\documentclass[a4paper,twoside,reqno,12pt]{amsart}

\usepackage{amsfonts,euscript,amsmath}
\usepackage{amssymb}

\def\AW*{{\sl AW*}-algebra}
\def\W*{{\sl W*}-algebra}
\def\C*{{\sl C*}-algebra}
\def\Cs*{{\sl C*}-subalgebra}
\def\AWm*{{\sl AW*}-module}
\def\AWsm*{{\sl AW*}-submodule}
\def\AF/{{\sl AF}-algebra}
\newcommand{\NN}{{\mathbb{N}}}
\newcommand{\RR}{{\mathbb{R}}}
\newcommand{\CC}{{\mathbb{C}}}
\newcommand{\Calgc}{{\sl C*}-algebra}
\newcommand{\Calg}{{\sl C*}-algebra }
\newcommand{\Calgs}{{\sl C*}-algebras }
\newcommand{\Calgsc}{{\sl C*}-algebras}
\newcommand{\Csalg}{{\sl C*}-subalgebra }

\newcommand{\Csalgs}{{\sl C*}-subalgebras }

\newcommand{\clos}[1]{{\kern.07em{}^c\kern-.11em{#1}}}  
\newcommand{\closa}{{\clos{\kern-.05em{A}}}}           
\newcommand{\rest}[2]{{{#1}_{\kern-.5pt|{#2}}}}  
\newcommand{\botimes}{{{}\,\mathop{\overline{\phantom{\bar v}}\kern-7.3pt\otimes}{}}}

\newcommand{\eps}{{\varepsilon}}
\renewcommand{\epsilon}{{\eps}}
\newcommand{\mloc}{{M_{\text{\rm loc}}(A)}}

\newcommand{\Mloc}[1]{{M_{\text{\rm loc}}({#1})}}
\newcommand{\MLOC}[1]{{M_{\text{\rm loc}}\bigl({#1}\bigr)}}

\newcommand{\qmax}{{Q_{\text{\rm max}}(A)}}
\newcommand{\Qmax}[1]{{Q_{\text{\rm max}}({#1})}}

\newcommand{\Qmaxs}[1]{{Q_{\text{\rm max}}^s({#1})}}
\newcommand{\qmaxs}{{\Qmaxs A}}
\newcommand{\Qmaxsb}[1]{{Q_{\text{\rm max}}^s({#1})_b}}
\newcommand{\qmaxsb}{{\Qmaxsb A}}

\newcommand{\CBsA}{C\!B_A}

\newcommand{\CCsAB}{{C\mkern-1mu C_{A\text{-}B}}}

\newcommand{\CCsCA}{{C\mkern-1mu C_{\CC\text{-}A}}}
\newcommand{\ela}{{{\mathcal E}\kern-2.1pt\ell(A)}}   
\newcommand{\longrightarrowraised}{{\hbox{\raise.5\jot\hbox{\scriptsize$\longrightarrow$}}}}
\newcommand{\longleftarrowraised}{{\hbox{\raise.5\jot\hbox{\scriptsize$\longleftarrow$}}}}
\newcommand{\dirlim}{{\smash{\underset{\longrightarrowraised}
                   {\operatorname{lim}}}\vphantom{a^{}_{a_f^{}}}}}
\newcommand{\smalllongrightarrowraised}{{\hbox{\hbox{\scriptsize$\rightarrow$}}}}
\newcommand{\smalldirlim}{{\smash{\underset{\smalllongrightarrowraised}
                   {\operatorname{lim}}}\vphantom{a^{}_{a_f^{}}}}}
\newcommand{\Dirlim}[1]{{\smash{\dirlim{}}\sp{}_{\,{#1}}
                   \vphantom{a^{}_{a_f^{}}}}}
\newcommand\alglim{{\smash{\underset{\longrightarrowraised}
                   {\operatorname{alg\,lim}}}\vphantom{a^{}_{a_f^{}}}}}
\newcommand\Alglim[1]{{\smash{\alglim{}}\sp{}_{\,{#1}}
                   \vphantom{a^{}_{a_f^{}}}}}
\newcommand{\invlim}{{\smash{\underset{\longleftarrowraised}
                   {\operatorname{lim}}}\vphantom{a^{}_{a_f^{}}}}}
\newcommand{\Invlim}[1]{{\smash{\invlim{}}\sp{}_{\,{#1}}
                   \vphantom{a^{}_{a_f^{}}}}}

\newcommand{\cat}[1]{{\mathcal{#1}}}
\newcommand{\Oone}{{\mathcal O}_1}
\newcommand{\AOoneB}{A\text{-}{\mathcal O}_1\text{-}B}
\newcommand{\AOoneC}{A\text{-}{\mathcal O}_1\text{-}\CC}
\newcommand{\COoneA}{\CC\text{-}{\mathcal O}_1\text{-}A}
\newcommand{\COA}{\CC\text{-}{\mathcal O}\text{-}A}
\newcommand{\Sone}{{\mathcal S}_1}
\newcommand{\ModR}{{\mathcal Mod}\text{-}R}
\newcommand{\BanA}{{\mathcal Ban}\text{-}A}

\newcommand{\Ir}{{{\mathfrak I}_r}}
\newcommand{\Icr}{{{\mathfrak I}_{cr}}}
\newcommand{\Ice}{{{\mathfrak I}_{ce}}}
\newcommand{\Ier}{{{\mathfrak I}_{er}}}
\newcommand{\Icer}{{{\mathfrak I}_{cer}}}
\newcommand{\Heress}{{{\mathfrak H}_e}}
\newcommand{\id}[1]{{\text{\rm id}_{{#1}}}}
\newcommand{\ol}[1]{{\overline{{#1}}}}
\newcommand{\Abar}{{\hbox{\kern.22em$\ol{\phantom{I}}\kern-.75em A$\kern.09em}}}
\newcommand{\smallAbar}{{\hbox{\scriptsize\kern.22em$\ol{\phantom{I}}\kern-.75em A$\kern.09em}}}
\newcommand{\Abarsa}{\Abar_{\!sa}}
\newcommand{\Ahat}{\hat A}
\newcommand{\Bbar}{{\hbox{\kern.22em$\ol{\phantom{I}}\kern-.75em B$\kern.09em}}}
\newcommand{\thnorm}{{|\mkern-2mu|\mkern-2mu|}}
\newcommand{\cbn}[1]{{\|{#1}\|_{cb}}}
\newcommand{\seqn}[1]{{({#1}_n)_{n\in\NN}}}
\newcommand{\LB}{{\mathcal L}_B}
\newcommand{\KB}{{\mathcal K}_B}
\def\sigfinite/{{$\sigma$-finite}}

\def\rightbrokenarrow{{\relbar\relbar\rightarrow}}
\def\longerrightarrow{{\setbox0=\hbox{$\rightbrokenarrow$}
                         \hbox to\wd0{\rightarrowfill}}}
\def\specialsixdiag#1#2#3#4#5#6#7#8{{\normalbaselines{\baselineskip20pt
                                 \lineskip3pt \lineskiplimit3pt}
     \def\mapfirstrighta{\smash{
                    \mathop{\longerrightarrow}\limits^{{#2}}}}
     \def\mapsecondrighta{\smash{
                    \mathop{\longerrightarrow}\limits^{{#4}}}}
     \def\mapfirstrightb{\smash{
                    \mathop{\longerrightarrow}\limits_{{#2}}}}
     \def\mapsecondrightb{\smash{
                    \mathop{\longerrightarrow}\limits_{{#4}}}}
     \def\mapdownl{\Big\downarrow\llap{$\vcenter{
                       \hbox{$\scriptstyle{#6\ \;}$}}$}}
     \def\mapdownm{\Big\downarrow\rlap{$\vcenter{
                       \hbox{$\scriptstyle{#7}$}}$}}
     \def\mapdownr{\Big\downarrow\rlap{$\vcenter{
                       \hbox{$\scriptstyle{#8}$}}$}}
     \begin{matrix}
            \vphantom{T_T}{#1}&\mapfirstrighta&{#3}&\mapsecondrighta&{#5}\\
            \mapdownl&{}&\mapdownm&{}&\mapdownr\\
            \vphantom{T^{T^T}}{#1}&\mapfirstrightb&{#3}&\mapsecondrightb&{#5}\\
     \end{matrix}}}
\def\specialsixdowndiag#1#2#3#4#5#6#7#8#9{{\normalbaselines{\baselineskip20pt
                                            \lineskip3pt \lineskiplimit3pt}
     \def\maprighta{\smash{
                    \mathop{\longerrightarrow}\limits^{{#2}}}}
     \def\maprightb{\smash{
                    \mathop{\longerrightarrow}\limits_{{}}}}
     \def\mapdownla{\Big\downarrow\llap{$\vcenter{
                       \hbox{$\scriptstyle{#4\ \;}$}}$}}
     \def\mapdownra{\Big\downarrow\rlap{$\vcenter{
                       \hbox{$\scriptstyle{#4}$}}$}}
     \def\mapdownlb{\Big\downarrow\llap{$\vcenter{
                       \hbox{$\scriptstyle{#6\ \;}$}}$}}
     \def\mapdownrb{\Big\downarrow\rlap{$\vcenter{
                       \hbox{$\scriptstyle{}$}}$}}
     \begin{matrix}
            \vphantom{\Big|}{#1}&\maprighta&{#3}&\mkern-10mu=\LB\bigl(C_b(U,H)\bigr)\\
            \mapdownla&{#5}&\mapdownra\\
            \vphantom{\Big|}{#7}&\maprightb&{#8}&\mkern-10mu=\mloc\hfill\\
            \mapdownlb&{#5}&\mapdownrb\\
            \vphantom{\Big|}{#7}&\maprightb&{#9}&\mkern-10mu=I(A)\hfill\\
     \end{matrix}}}
\theoremstyle{remark}
\newtheorem{definition}{\bf Definition}[section]
\newtheorem{nextdefinition}[definition]{\bf Definition}
\newtheorem{remark}[definition]{Remark}
\newtheorem{remarks}[definition]{Remarks}
\newtheorem{example}[definition]{Example}
\newtheorem{question}[definition]{Question}
\theoremstyle{plain}
\newtheorem{lemma}[definition]{Lemma}
\newtheorem{proposition}[definition]{Proposition}
\newtheorem{theorem}[definition]{Theorem}
\newtheorem{corollary}[definition]{Corollary}
\newenvironment{proofof[1]}{\em Proof of {#1}.\quad}{\quad\qed}
\newenvironment*{notation}{\medskip\noindent{\bf Notation.}}{\medskip}

\begin{document}

\title{Maximal C*-Algebras of Quotients and\\ Injective Envelopes of C*-Algebras}

\author{Pere Ara}
\address{Departament de Matem\`atiques, Universitat Aut\`onoma de Barcelona,
         E--08193 Bellaterra (Barcelona), Spain}
\email{para@mat.uab.cat}

\author{Martin Mathieu}
\address{Department of Pure Mathematics, Queen's University Belfast,
         Belfast BT7 1NN, Northern Ireland}
\email{m.m@qub.ac.uk}

\thanks{This paper is part of a research project supported by the Royal Society.
        The first-named author's research was partially supported by
        the DGI and the European Regional Fund, jointly, through
        project MTM2005-00934 and, in addition, by the Comissionat
        per Universitats i Recerca de la Generalitat de Catalunya.
        The second-named author gratefully acknowledges the support by the
        Centre de Recerca Matem\`atica, Barcelona.}


\begin{abstract}
A new {\sl C*}-enlargement of a \C* $A$ nested between the local multiplier
algebra $\mloc$ of $A$ and its injective envelope $I(A)$ is introduced.
Various aspects of this maximal \C* of quotients, $\qmax$, are studied,
notably in the setting of \AW*s. As a by-product we obtain a new example
of a type~I \C* $A$ such that $\Mloc\mloc\ne\mloc$.
\end{abstract}


\maketitle

\pagestyle{myheadings}
\markboth{\sc Pere Ara and Martin Mathieu}{\sc Maximal {\sl C*}-Algebras of Quotients}

\section{Introduction}\label{sect:intro}

The basic theory of the local multiplier algebra $\mloc$ of a \C*~$A$
was set out in~\cite{AMMzth}. It emerged that this concept is very useful,
sometimes crucial, in the structural study of operators between \C*s, especially
when investigating properties that are linked to the ideal structure of~$A$.
In this sense, $\mloc$ plays the same role for non-simple \C*s as the multiplier
algebra $M(A)$ does for simple~$A$. For an example, we mention Sakai's theorem
stating that every derivation of a simple \C* $A$ becomes inner in $M(A)$ and
Pedersen's theorem stating that every derivation of a separable \C* $A$ becomes
inner in~$\mloc$; see \cite[Corollary~4.2.9]{AMMzth} and \cite[Theorem~4.2.20]{AMMzth},
respectively.

Of course, there are other enlargements of a \C*~$A$ relevant to operator theory on~$A$,
notably the injective envelope $I(A)$, which plays an important role for completely
bounded and completely positive operators. The relation between $\mloc$ and $I(A)$
was first studied in~\cite{Frank}, \cite{FrankPaulsen}, where it was shown that,
for every $k\geq1$,
the $k$-times iterated local multiplier algebra $M_{\text{\rm loc}}^{(k)}(A)$ of $A$
is contained in~$I(A)$; see also Proposition~\ref{prop:mloc-into-qmax} and
Theorem~\ref{thm:qmax-into-injenv} below. When $A$ is commutative, $\mloc$ is a commutative
\AW*, hence injective; as a result, $\mloc=I(A)$. In the non-commutative setting,
$\mloc$ is typically much smaller than $I(A)$ though it may happen that
$\Mloc\mloc=I(A)$; see, e.g., page~\pageref{p:mlocmloc} below.
In this paper, we shall introduce another \C* nested between $\mloc$ and $I(A)$,
which we call the maximal \C* of quotients of~$A$, and study some of its basic properties.

In analogy to the notion of a two-sided ring of quotients in non-commutative ring
theory, a \textit{\C* of quotients\/} of a given \C* $A$ can be thought of as a unital
\C* $B$, containing a *-isomorphic image of~$A$, such that
\vskip2pt
\begin{enumerate}
\item[(1)]$\{b\in B\mid bJ+b^*J\subseteq A\,\text{ for some }\,J\in{\mathfrak I}\}$
           is dense in~$B$,
\smallskip
\item[(2)]for all $b\in B$, $J\in{\mathfrak I}$:\quad$bJ=0\;{}\Longrightarrow{}\; b=0$,
\end{enumerate}
\vskip2pt\noindent
where ${\mathfrak I}$ is a `good' set of closed right ideals in~$A$.
A generally accepted notion of a `good' set of right ideals requires $\mathfrak I$
to be a \textit{Gabriel filter\/}, that is, a non-empty subset of $\Ir$, the lattice
of all right ideals of $A$, satisfying the following properties:
\vskip2pt
\begin{enumerate}
\item[\rm (i)] if $J_1,\,J_2\in{\mathfrak I}\,$ then $J_1\cap J_2\in{\mathfrak I}$;
\smallskip
\item[\rm (ii)] if $J\in{\mathfrak I}$, $I\in\Ir$ and $J\subseteq I$
                then $I\in{\mathfrak I}$;
\smallskip
\item[\rm (iii)] if $J\in{\mathfrak I}$, $x\in A$ then
                 $x^{-1}J=\{a\in A\mid xa\in J\}\in{\mathfrak I}$;
\smallskip
\item[\rm (iv)] if $J\in{\mathfrak I}$, $I\in\Ir$ and $x^{-1}I\in{\mathfrak I}\,$
                for all $x\in J$ then $I\in{\mathfrak I}$.
\end{enumerate}
\vskip2pt\noindent
For more details we refer to~\cite{Sten}.
For instance, $\Ier$ consisting of all essential right ideals is a Gabriel filter.
(Here, a right ideal $J$ of $A$ is said to be \textit{essential\/} if it intersects
non-trivially with every other non-zero right ideal of~$A$.)

Replacing $\Ir$ by $\Icr$, the lattice of all \textit{closed\/} right ideals of $A$, and
restricting our attention to Gabriel filters on~$\Icr$
we can choose for ${\mathfrak I}$ the trivial filter ${\mathfrak I}=\{A\}$;
in this case we get the multiplier algebra $M(A)$ as a \C* of quotients.
With respect to ${\mathfrak I}=\Icer$, the filter of all
closed essential right ideals of $A$, the new \C* introduced in this paper, $\qmax$,
turns out to be the maximal \C* of quotients.
However, essentiality of a \textit{closed\/} right ideal could well be defined in
various other ways depending on the category one has in mind, and the subtle relations
between the (at least) five possible interpretations will occupy us throughout this work.

It is fairly evident that $\mloc$ can be embedded into $\qmax$
(Proposition~\ref{prop:mloc-into-qmax}), and it is easy to
give examples where the two \C*s are different. Let $A$ be a unital simple separable
and infinite-dimensional \C*. Then $\mloc=A$. Let $J$ be a maximal right ideal of~$A$.
Then $J$ is essential, since $A$ does not contain minimal right ideals by assumption.
As $\qmax$ contains the multiplier algebra of every essential hereditary \Cs* of $A$
(see Remark~\ref{rem:ess-right-ideals-vs-algebras}), in particular, $M(J\cap J^*)$
and the latter is non-separable, $\qmax\ne A$ in this case.

We shall introduce $\qmax$ in Section~\ref{sect:max-quotients} in analogy
to the construction of the maximal symmetric ring of quotients, which already came up
in work by Utumi and Schelter but was studied in more detail recently by
Lanning~\cite{Lann}. In doing this we shall follow the algebraic approach
that was pursued to construct the local multiplier algebra in~\cite{AMMzth}
and which has many advantages over the direct limit construction in the category
of \C*s. In fact, it turns out that it is impossible to give an analogous
direct limit description for $\qmax$. The basic idea is to introduce a suitable
order structure on $\qmaxs$, the maximal symmetric algebra of quotients of~$A$,
and to define order-bounded elements which, when collected together, form a
pre-\C* under the order-unit norm the completion of which is $\qmax$. This
approach was already outlined in~\cite{Arazth}.

Focussing on the interrelation between $\qmax$ and the injective envelope $I(A)$
in the remainder of Section~\ref{sect:max-quotients} we find that $\qmax$ is
contained in the regular monotone completion of $A$ in the sense of
Hamana~\cite{Hameoa} in Proposition~\ref{prop:qmax_into_abar} and thus that
$\qmax=A$, whenever $A$ is monotone complete. We also obtain the maximal \C* of quotients
as a direct limit in the category of operator modules (Corollary~\ref{cor:eb-into-injenv}).

Because of its importance for our approach we will review the concept of the injective
envelope and its basic properties in Section~\ref{sect:inj-envel}. We will follow a
categorical approach and work entirely in the category $\Oone$ of operator spaces
and complete contractions. This category also plays an important role in
Section~\ref{sect:ess-right-ideals}, where we discuss five possibly different
notions of essentiality for one-sided closed ideals in a \C*.

It is not too difficult to verify that $\mloc=A$ for every \AW*~$A$; see
\cite[Theorem~2.3.8]{AMMzth}. Since it is unknown whether every \AW* $A$ is monotone
complete, the corresponding question for $\qmax$ remains open.
In  Section~\ref{sect:awstars}, we give a positive answer in the finite case
(Theorem~\ref{thm:qmax-finite-aw}) and under the assumption that every orthogonal
family of projections is at most countable (Corollary~\ref{cor:qmax-sigf-aw}).
We also prove that such a $\sigma$-finite \AW* $A$ is $A$-injective and that
every closed right ideal of $A$ with zero left annihilator is $\Oone$-essential
(Proposition~\ref{prop:sigf-aw-selfinj} and
Corollary~\ref{cor:ess-right-ideals-in-sigfinite-aw}, respectively).
In Proposition~\ref{prop:qmax_vN} we obtain that, for a separable \C* $A$,
$\qmax$ is a von Neumann algebra if and only if $A$ contains a minimal essential
ideal consisting of compact elements, an analogue to the recent result by
Argerami and Farenick for $\mloc$ and~$I(A)$~\cite{ArgFar1}.

The example of $A=C(X)\otimes B(H)$, where $X$ is a compact Hausdorff space and
$H$ is a separable Hilbert space, will be studied in detail in Section~\ref{sect:example}.
We determine $\mloc$ in this case, and also $\qmax$ and $I(A)$ under the additional
hypothesis that $X$ is Stonean. Not surprisingly the notion of {\sl AW*}-modules
and various conditions on the topology of $X$ will play a vital role.

In his seminal paper~\cite{Pedsve}, where he introduced $\mloc$ as a direct limit of
multiplier algebras of essential closed ideals of~$A$, Pedersen asked whether there is
a \C* $A$ with the property that $\mloc\ne\Mloc\mloc$. The first example of such a \C*
was given in~\cite{AMMzs}; in fact, a class of unital prime \AF/s was exhibited.
Note that each example of this type has to be antiliminal for the following reason.
Let $A$ be a prime \C* with the property that its largest postliminal ideal $I$ is non-zero.
As $\mloc=\Mloc I$ we can replace $A$ by $I$,
if necessary, and thus assume that $A$ is type~I. It is well known that $A$ has to
contain a copy of the compact operators as an ideal in this case.
Hence $\mloc$ is isomorphic to $B(H)$ and thus $\Mloc\mloc=\mloc$.
At the end of Section~\ref{sect:example} we shall discuss a new example of a type~I
\C* $A$ such that $\Mloc\mloc\ne\mloc$ (Corollary~\ref{cor:type-one-mloc}).
A related example was found independently in~\cite{ArgFar2} by other methods;
one of the differences between the two examples is that the centre of ours has a Stonean
structure space.

The paper concludes with a list of open problems that illustrate why the study
of even basic properties of the maximal \C* of quotients appears to be somewhat
more involved than that of the local multiplier algebra.

\section{Injective Envelopes}\label{sect:inj-envel}

We begin this section by reviewing the concept of injective envelope of a \C*,
and indeed of an operator space, from a categorical point of view.

Let $\cat C$ be a category, and let $\cat H$ be a class of morphisms in~$\cat C$
which is closed under composition and contains all isomorphisms of~$\cat C$.
We shall denote by $\cat C(A,B)$ the set of all morphisms from an object $A$ to an object
$B$ in~$\cat C$.

\begin{definition}
An object $I$ in $\cat C$ is called {\em$\cat H$-injective\/} if the mapping
\[
\cat C(h,I)\colon\cat C(B,I)\to\cat C(A,I),\quad g\mapsto gh
\]
is surjective for all morphisms $h\in\cat C(A,B)$ in~$\cat H$.
A morphism $h\in\cat H$ is called {\em$\cat H$-essential\/} if, for every morphism
$g$, $gh\in\cat H$ if and only if $g\in\cat H$. We will denote the class of all
$\cat H$-essential morphisms by~$\cat H^*$.

Given an object $A$ in $\cat C$, a {\em$\cat H$-injective envelope\/} of $A$
is a $\cat H$-injective object $I(A)$ together with a $\cat H$-essential morphism
$\iota_A\colon A\to I(A)$.
\end{definition}

One expects injective envelopes to be unique up to isomorphism; indeed, a stronger
property is true.

\begin{proposition}
Every morphism in\/ $\cat H$ between two\/ $\cat H$-injective envelopes of the same object
is an isomorphism.
\end{proposition}

\begin{proof}
Let $h\colon A\to I$ and $k\colon A\to J$ be two $\cat H$-injective envelopes
of the object $A$ in~$\cat C$. Since $J$ is $\cat H$-injective, $k$ factors as
$k=fh$, where $f\in\cat H$ because $k\in\cat H$ and $h\in\cat H^*$.
Since $I$ is $\cat H$-injective, $\text{id}_I$ factors as $\text{id}_I=gf$,
wherefore $h=gfh=gk$. As before, this entails that $g\in\cat H$.
Applying the $\cat H$-injectivity of $J$ once again, we can factor $\text{id}_J$
as $\text{id}_J=lg$. As a result, $g$ is an isomorphism between $J$ and~$I$.

Clearly, the second half of the above argument gives the statement for an arbitrary
$f\in\cat C(I,J)$ which belongs to~$\cat H$.
\end{proof}

Note, however, that $\cat H$-injective envelopes may not be unique up to {\it unique\/}
isomorphism in general, and may not be natural~\cite{AHRT}.

In the construction of injective envelopes of operator spaces, the following notions
will be helpful.

\begin{nextdefinition}
Let $\iota\in\cat C(A,I)$; we say that $\iota$ is {\it rigid\/} if $\id I$ is
the unique endomorphism extending~$\id A$. That is, whenever $g\in\cat C(I,I)$
is such that $g\iota=\iota$ then $g=\id I$. We say that $\cat C$ admits
{\em rigid $\cat H$-injective envelopes\/} if, for every object $A$ in~$\cat C$,
there is a $\cat H$-injective envelope $I(A)$ such that $\iota_A\colon A\to I(A)$
is rigid.
\end{nextdefinition}

We note that, if $\cat C$ admits rigid $\cat H$-injective envelopes, then
$\cat H$-injective envelopes are unique up to unique isomorphism. Indeed,
given an isomorphism $g\colon A_1\to A_2$ there is a unique isomorphism
$\tilde g\colon I(A_1)\to I(A_2)$ such that $\tilde g\iota_{A_1}=\iota_{A_2}g$.

\medskip
In order to apply the above discussion to operator spaces and operator systems,
we first fix our terminology.

\begin{notation}
We denote by $\Oone$ the category of operator spaces with complete contractions
as morphisms; in this category, $\cat H$ consists of the complete isometries.
We denote by $\Sone$ the category of operator systems with completely positive
unital linear maps; in this category, $\cat H$ consists of the unital complete
isometries.
\end{notation}

Before we turn our attention to the existence of injective envelopes, we first note
that $\cat H$-injective objects in the subcategory $\Sone$ are the same as in the
category~$\Oone$.

\begin{proposition}
The following conditions on an operator system\/ $S\subseteq B(H)$ are equivalent:
\begin{enumerate}
\item[{\rm(a) }] $S$ is injective in\/ $\Oone$;
\item[{\rm(b) }] $S$ is injective in\/ $\Sone$;
\item[{\rm(c) }] there is a completely positive projection\/ $\varphi$ from\/ $B(H)$
                 onto\/~$S$.
\end{enumerate}
\end{proposition}

\begin{proof}
Since the completely positive unital maps are precisely the unital complete contractions
\cite[Propositions~3.5 and~3.6]{Paul}, the implication (a)${}\Rightarrow{}$(b) is clear.
Assuming (b) we get a completely positive unital extension $\varphi$ of $\id S$ to $B(H)$,
which is thus a projection onto~$S$. Finally, assuming~(c), let $f\in\Oone(E,S)$ for some
operator space~$E$. If $h\colon E\to F$ is a complete isometry into the operator space
$F$, by Wittstock's extension theorem \cite[Theorem~8.2]{Paul}, there is a complete
contraction $\tilde g\colon F\to B(H)$ such that $\tilde gh=j_Sf$, where $j_S\colon S\to B(H)$
is the canonical monomorphism. Consequently, $g=\varphi\tilde g$ yields the desired
extension of $f$ to $F$ into~$S$.
\end{proof}

In the following, we will therefore unambiguously call $\cat H$-injective envelopes
simply {\em injective envelopes\/} in both $\Oone$ and~$\Sone$. Note, however, that
the injective envelope of an operator system is an operator system whereas the injective
envelope of an operator space is an operator space.

The \Calg $B(H)$ is injective in $\Oone$ (by Wittstock's theorem) and in $\Sone$ (by
Arveson's theorem \cite[Theorem~7.5]{Paul}). The idea of an injective envelope of
an operator space $E$ is to find a ``smallest'' injective object in $\Oone$ containing~$E$.
This goes back to work by Hamana~\cite{Hamsvn}; our approach is different from his
and also Paulsen's in \cite{Paul} and allows us to shorten the exposition slightly.

We start with the following observation.

\begin{lemma}\label{lem:rigid1}
Let\/ $\iota\colon E\to I$ be a rigid completely isometric mapping from the
operator space\/ $E$ into the injective operator space~\/$I$. Then\/ $\iota$ is essential
and hence\/ $I$ is an injective envelope of~\/$E$.
\end{lemma}

\begin{proof}
Let $g\colon I\to F$ be a morphism in $\Oone$ such that $g\iota$ is a complete
isometry. Since $I$ is injective, there is a morphism $h\colon F\to I$ such that
$hg\iota=\iota$. This implies that $hg$ is an extension of $\id E$, wherefore the rigidity
of $\iota$ entails that $hg=\id I$. But this implies that $g$ is a complete isometry,
whence $\iota$ is essential.
\end{proof}

\begin{theorem}\label{thm:injenv-exist}
For every operator space\/ $E$ there exists a rigid injective envelope\/ $I(E)$;
this is uniquely determined up to unique completely isometric isomorphism.
\end{theorem}

\begin{proof}
By Ruan's theorem \cite[Theorem~13.4]{Paul}, there exists a complete isometry $j$
from $E$ onto an operator subspace $F$ of some~$B(H)$. The first two paragraphs
of the proof of Theorem~15.4 in~\cite{Paul} (on page 210) give the existence of a minimal
$F$-projection $\varphi\colon B(H)\to B(H)$; since $\varphi$ is completely contractive
and $B(H)$ is injective, it follows that $I=\varphi(B(H))$ is an injective containing~$F$.
Lemma~15.5 in~\cite{Paul} states that $\rest{\varphi}{F}\colon F\to I$
is rigid. Therefore the complete isometry $\iota=\rest{\varphi}{F}\circ j\colon E\to I$
is rigid and  Lemma~\ref{lem:rigid1} entails that $I$ is an injective envelope of~$E$.

We have already remarked that the existence of rigid injective envelopes makes them
unique up to unique isomorphism; thus we will henceforth denote this unique rigid
injective envelope of $E$ by~$I(E)$.
\end{proof}

We next note that a rigid morphism into an injective object can be factorised only
into rigid components.

\begin{lemma}\label{lem:rigid2}
Let\/ $\iota\colon E\to I$ be a rigid morphism in\/ $\Oone$ into an injective
operator space. If\/ $\iota=kh$ for two complete isometries\/ $k$ and~\/$h$, then
both\/ $k$ and\/ $h$ have to be rigid.
\end{lemma}

\begin{proof}
Suppose that $h\in\cat H(E,F)$ and $k\in\cat H(F,I)$, that is, both are complete
isometries, are given with the property $\iota=kh$. Let $f$ be an endomorphism of $I$
such that $fk=k$. Then
\[
f\iota=fkh=kh=\iota
\]
so that the rigidity of $\iota$ entails that $f=\id I$. Hence $k$ is rigid.

Suppose that $g$ is an endomorphism of $F$ such that $gh=h$. Then we have the
following commutative diagram
\medskip
\[
\specialsixdiag{E}{h}{F}{k}{I}{\id A}{g}{\tilde g}
\]

\medskip
\noindent
where the injectivity of $I$ yields the endomorphism $\tilde g$ with the property
$\tilde gk=kg$, since $k\in\cat H$. Consequently,
\[
\tilde g\iota=\tilde gkh=kgh=kh=\iota
\]
so that rigidity of $\iota$ entails that $\tilde g=\id I$. Thus $k=kg$ and $g=\id F$,
as $k$ is injective. Hence $h$ is rigid.
\end{proof}

We now introduce two notions that are, in a way, dual to each other.

\begin{nextdefinition}
Let $I$ be an injective object in $\Oone$. Let $\iota\in\cat H(E,I)$.
Then $I$ is a {\em minimal injective containing\/} $E$ if, whenever $\iota=kh$ with
$h\in\cat H(E,I_1)$ for some injective $I_1$ in $\Oone$ and $k\in\cat H(I_1,I)$,
then $k$ is an isomorphism.
\end{nextdefinition}

\begin{nextdefinition}
Let $I$ be an object in $\Oone$. Let $\iota\in\cat H^*(E,I)$.
Then $I$ is a {\em maximal essential for\/} $E$ if, whenever
$h\in\cat H(E,F)$ is essential, then there exists $k\in\cat H(F,I)$ such that
$\iota=kh$.
\end{nextdefinition}

As we shall see now these properties characterise (rigid) injective envelopes.

\begin{theorem}\label{thm:injenv-char}
Let\/ $E$ and\/ $I$ be operator spaces, and let\/ $\iota\colon E\to I$
be a complete isometry. The following conditions are equivalent:
\begin{enumerate}
\item[{\rm(a) }] $\iota$ is rigid and\/ $I$ is injective;
\item[{\rm(b) }] $I$ is an injective envelope of\/ $E$;
\item[{\rm(c) }] $I$ is a minimal injective containing\/ $E$;
\item[{\rm(d) }] $I$ is a maximal essential for\/ $E$;
\item[{\rm(e) }] $I$ is a maximal essential injective for\/ $E$.
\end{enumerate}
\end{theorem}

\begin{proof}
(a)${}\Rightarrow{}$(b)\enspace
This is Lemma~\ref{lem:rigid1}.

\noindent
(b)${}\Rightarrow{}$(a)\enspace
The injective envelope $I$ is (uniquely) completely isometrically isomorphic to the
rigid injective envelope $I(E)$ constructed in Theorem~\ref{thm:injenv-exist}; hence
it is rigid itself.

\noindent
(a)${}\Rightarrow{}$(c)\enspace
Suppose that $I_1$ is injective and that $h\colon E\to I_1$ and $k\colon I_1\to I$
are complete isometries such that $\iota=kh$. By Lemma~\ref{lem:rigid2}, $h$ is rigid
and thus essential, by Lemma~\ref{lem:rigid1}. Therefore $I_1$ is an injective
envelope of $E$ and the uniqueness of $I(E)$ yields that $k$ is an isomorphism.
Hence, $I$ is minimal injective.

\noindent
(c)${}\Rightarrow{}$(d)\enspace
Suppose $h\colon E\to F$ is a complete isometry. Then there is a
complete contraction $k\colon F\to I$ such that $kh=\iota$. If $h$ is essential,
then $k$ is completely isometric, wherefore $I$ is maximal essential.

\noindent
(d)${}\Rightarrow{}$(e)\enspace
Suppose $\iota\colon E\to I$ is maximal essential; then, by definition, there is
a complete isometry $k\colon I(E)\to I$ such that $k\iota_E=\iota$. On the other
hand, by injectivity, there is a complete isometry $h\colon I\to I(E)$ such that
$h\iota=\iota_E$. Since $\iota_E$ is rigid, it follows that $hk=\id{I(E)}$.
This entails that $h(\id I-kh)=0$ which yields $kh=\id I$, as $h$ is injective.
Consequently, $I$ and $I(E)$ are completely isometrically isomorphic wherefore
$I$ is injective.

\noindent
(e)${}\Rightarrow{}$(b)\enspace
Since $I$ is essential and injective, it must be the uniquely determined injective
envelope of~$E$.
\end{proof}

\begin{corollary}\label{cor:injenv-c*}
Every operator system\/ $S$ has an up to unique completely isometric isomorphism uniquely
determined injective envelope\/ $I(S)$, which is an operator system, and for every unital
\Calg $A$ we can choose an injective envelope\/ $I(A)$ which is a unital \C* and contains\/
$A$ as a unital \C*.
\end{corollary}

\begin{proof}
The injective envelope of $S$ in $\Oone$ constructed in Theorem~\ref{thm:injenv-exist}
is completely isometric to $\varphi(B(H))$, which is an operator system.
If $S=A$ is a unital \C*, then $\varphi(B(H))$ can be given the structure of a
unital \C* via the product $a\circ b=\varphi(ab)$ (where we, for simplicity,
assume $A$ embedded into $B(H)$), by the Choi--Effros theorem \cite[Theorem~15.2]{Paul}.
As $\varphi$ is the identity on $A$, $A$ is a unital \Cs* of~$I(A)$.
\end{proof}

Rigidity can be exploited to obtain properties of the injective envelope readily, as in
the following two results due to Hamana~\cite{Hamsvn}.

\begin{proposition}
Let\/ $J$ be a non-zero closed ideal in\/ $I(A)$ for some unital \C*~$A$.
Then\/ $J\cap A\ne0$.
\end{proposition}

\begin{proof}
Let $\pi\colon I(A)\to I(A)/J$ denote the canonical epimorphism. If $J\cap A=0$,
then $\pi\iota_A$ is a *-monomorphism, hence completely isometric. The injectivity
of $I(A)$ thus yields a complete contraction $\tau\colon I(A)/J\to I(A)$ such that
$\iota_A=\tau\pi\iota_A$. Rigidity of $\iota_A$ entails that $\tau\pi=\id{I(A)}$,
wherefore $J=0$.
\end{proof}

\begin{corollary}
Let\/ $A$ be a simple unital \C* (a prime \C*, respectively).
Then\/ $I(A)$ is simple (prime, respectively).
\end{corollary}

Note, however, that $K(H)$ with $H$ infinite dimensional provides an example of
a non-unital simple \C* such that its injective envelope, which is $B(H)$, is not simple.

\smallskip

Inside the injective envelope of a \C* $A$ we find an important \Cs* containing~$A$
that has a particularly well-behaved order structure, the regular monotone completion of~$A$,
which we will now review.

Recall that a \C* $A$ is said to be \textit{monotone complete\/} if every bounded increasing
net in the selfadjoint part $A_{sa}$ has a supremum. A subset $S$ of $A_{sa}$ in a
monotone complete
\C* $A$ is called \textit{monotone closed\/} if the supremum of every bounded increasing net
in $S$ and the infimum of every bounded decreasing net in $S$ belong to~$S$. The
\textit{monotone closure\/} of a \Cs* $B$ in $A$ is given by
$C+iC$, where $C$ is the smallest monotone closed subset of $A_{sa}$ containing~$B_{sa}$.
(It is readily seen that $C$ is a real linear subspace of~$A$.)
The \C* $B$ is said to be \textit{monotone closed\/} (\textit{in\/} $A$)
if it agrees with its monotone closure.

The injective envelope $I(A)$ of a unital \C* $A$ is monotone complete.
Indeed, in the construction of Corollary~\ref{cor:injenv-c*} every bounded
increasing net in $I(A)_{sa}$ has a supremum in $B(H)_{sa}$. The completely
positive mapping $\varphi\colon B(H)\to I(A)$ maps this supremum to the
supremum in~$I(A)_{sa}$.
The \textit{regular monotone completion\/} $\Abar$ of $A$ is defined as the monotone
closure of $A$ inside~$I(A)$. It is characterised by the following properties.
\begin{enumerate}
\item $\Abar$ is monotone complete;
\item $A$ is order dense in $\Abar$, that is, for each $x\in\Abarsa$ we have
      \[x=\sup\{a\in A_{sa}\mid a\leq x\};\]
\item the inclusion $A\subseteq\Abar$ is sup-preserving, that is, if $S\subseteq A_{sa}$
      has a supremum $s\in A_{sa}$ then $s=\Abarsa\text{-}\sup S$.
\end{enumerate}
See \cite[Theorem~3.1]{Hameoa} for details. If $A$ is separable then $\Abar=\Ahat$, the
regular $\sigma$-completion of $A$ in the sense of Wright~\cite{Wrighteza}.
Without any restriction on~$A$ we have $\ol{M_n(A)}=M_n(\Abar)$ for all $n\in\NN$
\cite[Corollary~3.10]{Hameoa} and $\ol{pAp}=p\Abar p$ for every projection
$p\in A$ \cite[Proposition~1.11]{Hameoa}.

We shall frequently use the following proposition, often without making
explicit mention of it.

\begin{proposition}\label{prop:invariant-injective}
Let\/ $A$ be a unital \C* and let\/ $B$ be \Cs* of\/ $I(A)$ containing~$A$.
Then\/ $I(A)=I(B)$. If, moreover, $B\subseteq\Abar$ then\/ $\Abar=\Bbar$.
\end{proposition}
\begin{proof}
The inclusion $B\hookrightarrow I(A)$ is rigid in $\Oone$, since $A\subseteq B$
and $I(A)$ is a rigid injective envelope of~$A$. Lemma~\ref{lem:rigid1} thus entails
that $I(A)$ is an injective envelope of~$B$, hence $I(A)=I(B)$. Suppose that
$B\subseteq\Abar$; then certainly $\Bbar\subseteq\Abar$.
Since $A\subseteq B$, the converse inclusion $\Abar\subseteq\Bbar$ holds as well.
\end{proof}
Hamana shows in \cite[Theorem~6.6]{Hameoa} that, if $A$ is a type~I \C*, then
$\Abar=I(A)$ is a type~I \AW*.

An important concept Hamana introduces in~\cite{Hameta} is the notion of an open projection
relative to $\Abar$. A projection $p\in\Abar$ is called \textit{open\/} if $p$ is the
supremum in $\Abarsa$ of an increasing net of positive elements in~$A$. The complementary
projection $1-p$ is called \textit{closed\/} in this case. The relation with the `usual'
kind of open projections in the second dual $A^{**}$ of $A$ in the sense of Akemann
and Giles--Kummer, see \cite[3.11.10]{Pedsvn}, is as follows. Consider $A$ as a \Cs*
of $B(H_u)$, where $H_u$ is the Hilbert space of the universal representation.
Then $A^{**}$ can be identified with the ultraweak closure of $A$ in~$B(H_u)$.
Construct $I(A)$ inside of $B(H_u)$ via a minimal $A$-projection~$\varphi$, see
Theorem~\ref{thm:injenv-exist} and Corollary~\ref{cor:injenv-c*}. Then the open projections
in $\Abar$ in Hamana's sense are the images of the open projections in $A^{**}$
under~$\varphi$; see \cite[Lemma~6.4]{Hametb}.

The open projections in $\Abar$ will play an important role in the next section,
especially because of their relation with one-sided ideals (cf.~Lemma~\ref{lem:openproj}).
For a hereditary \Csalg $B$ of $A$, there exists a unique open projection
$p_B\in\Abar$ with the property that, together with the canonical inclusions,
$p_BI(A)p_B$ is the injective envelope of~$B$ and $p_B\Abar p_B$ is the regular monotone
completion of~$B$, respectively. Moreover, $p_B=I(A)_{sa}\text{-}\sup e_\alpha$ for each
approximate identity $(e_\alpha)$ of~$B$ and $p_B$ is central if $B$ is a two-sided
ideal.\label{p:B}
See \cite[Lemma~1.1]{Hameta} and \cite[Theorem~6.5]{Hametb}.

\section{Essentiality of One-sided Ideals}\label{sect:ess-right-ideals}

There exist a number of possible concepts of essentiality for one-sided ideals
of \Calgsc. The relations between these seem to be not fully understood; we therefore
devote this section to a closer look at the various notions. Throughout, we will work
with right ideals; it is clear that an analogous discussion can be undertaken
for left ideals. Note that for two-sided closed ideals all these concepts agree with each other.

\begin{definition}
An operator subspace $E$ of an operator space $F$ is said to be {\it essential\/}
({\it in\/} $F$) provided the inclusion $E\hookrightarrow F$ is essential in~$\Oone$.
This is equivalent to the requirement that the map
\[
E\hookrightarrow F\overset{\iota_F}{\rightarrow}I(F)
\]
yields an injective envelope of~$E$.
\end{definition}

\begin{definition}
Let $R$ be a ring, and let $N\subseteq M$ be right $R$-modules. We say that
$N$ is {\it algebraically essential in\/} $M$ provided $N\cap K\ne0$ for
every non-zero submodule $K$ of~$M$. This is equivalent to the requirement
that $N\hookrightarrow M$ is essential in the category $\ModR$ of right $R$-modules
(where $\mathcal H$ of course consists of the injective $R$-module maps).
\end{definition}

For a subset $X$ of an algebra $A$ we denote by $\ell(X)$ and $r(X)$ its
{\it left\/} and its {\it right annihilator}, respectively; that is,
\[
\ell(X)=\{a\in A\mid aX=0\}\quad\text{and}\quad r(X)=\{a\in A\mid Xa=0\}.
\]
The {\it two-sided annihilator\/} of $X$ is defined by
\[
X^\perp=\{a\in A\mid aX=Xa=0\}.
\]
We shall emphasize the surrounding algebra in cases of ambiguity by writing
$\ell_A(X)$ etc.\label{p:annis}

For a two-sided closed ideal $I$ of a \C* $A$ we have $I^\perp=\ell(I)=r(I)$.

\begin{definition}
Let $J$ be a closed right ideal in a \Calg $A$. We say that $J$ is {\it KP-essential\/}
({\it essential in the sense of Kaneda--Paulsen~\cite{KanPaul}\/})
provided that $r(J)=0$. This is equivalent to the requirement that the two-sided
ideal $J^*J$ is essential as a two-sided ideal.
\end{definition}

As we shall see the last of the three notions of essentiality is the weakest one.
Indeed, if $J$ is a closed right ideal of the \C* $A$ then

\medskip
\begin{center}
$J$ algebraically essential
${}\Rightarrow{}$
$J$ $\Oone$-essential
${}\Rightarrow{}$
$\ell_{A}(J)=0$
${}\Rightarrow{}$
$J$ KP-essential.
\end{center}

\medskip
In order to discuss this in a comprehensive setting, we first introduce some
terminology and notation.
By an {\it operator $A$-$B$-bimodule\/} $E$ we understand an operator space $E$
which is a bimodule over the \Calgs $A$ and $B$ such that the module operations
are completely contractive. We shall denote the category of operator $A$-$B$-bimodules
together with completely contractive $A$-$B$-bimodule maps by $\AOoneB$ and the sets
of morphisms therein by $\CCsAB(E,F)$.
Objects in $\AOoneC$ and in $\COoneA$, respectively are called
{\it operator left $A$-modules\/} and {\it operator right $A$-modules}, respectively.
Clearly, $\CC\text{-}\Oone\text{-}\CC=\Oone$. Whenever the \Calgs are unital, we will
assume that the modules are unitary.

Part of the following result is contained in~\cite[Theorem~5.1]{KanPaul}.

\begin{theorem}\label{thm:ess-right-ideals}
Let\/ $J$ be a closed right ideal of a \C* $A$, and let\/ $(e_\alpha)$
be an approximate identity of\/ $J^*\cap J$. The following conditions are
equivalent:
\begin{enumerate}
\item[{\rm(a) }] $J$ is\/ $\Oone$-essential;
\vskip2pt
\item[{\rm(b) }] $J$ is\/ $\COoneA$-essential;
\vskip2pt
\item[{\rm(c) }] $\|(a_{kl})\|=\sup_\alpha\|(e_\alpha a_{kl})\|$\
                 for all\/ $(a_{kl})\in M_n(A)$, $n\in\NN$;
                 \vskip2pt
\item[{\rm(d) }] $\ell_{I(A)}(J)=0$.
\end{enumerate}
If\/ $A$ is unital then the following condition is equivalent as well:
\begin{enumerate}
\item[{\rm(e) }] $\CCsCA(A/J,I(A))=0$.
\end{enumerate}
\end{theorem}

\begin{proof}
\noindent
(a)${}\Rightarrow{}$(b)\enspace
This is clear.

\noindent
(b)${}\Rightarrow{}$(c)\enspace
Let us define $p\colon A\to\RR_+$ by $p(a)=\sup_\alpha\|e_\alpha a\|$, $a\in A$.
Then $p$ is a semi-norm such that $\ker p$ is a closed right $A$-module. Thus,
$\tilde A=A/\ker p$ becomes an operator $\CC$-$A$-bimodule with the matrix norms
given by
\begin{equation*}
\thnorm(a_{kl}+\ker p)\thnorm=\sup_\alpha\|(e_\alpha a_{kl})\|
\qquad\bigl((a_{kl})\in M_n(A),\,n\in\NN\bigr).
\end{equation*}
The canonical quotient map $\phi\colon A\to\tilde A$ is a completely contractive
$\CC$-$A$-bimodule mapping which is completely isometric on~$J$. Hence, by assumption,
$\phi$ is completely isometric on~$A$ which implies that
$\|(a_{kl})\|=\thnorm(a_{kl})\thnorm=\sup_\alpha\|(e_\alpha a_{kl})\|$
as claimed.

\noindent
(c)${}\Rightarrow{}$(d)\enspace
Set $F=\ell_{I(A)}(J)$; this is a closed operator subspace of~$I(A)$.
Considering the sequence
$\,A\overset{\iota_A}{\longrightarrow}I(A)\overset{\pi}{\longrightarrow}I(A)/F\,$
we have, by hypothesis,
\begin{equation*}
\|a+x\|\geq\sup_\alpha\|(a+x)e_\alpha\|=\sup_\alpha\|ae_\alpha\|=\|a\|
\qquad(a\in A,\ x\in F)
\end{equation*}
thus $\pi\circ\iota_A$ is an isometry. Using the hypothesis (c) and the same argument
on matrices, we find that $\pi\circ\iota_A$ is a complete isometry. Since $\iota_A$
is essential, we conclude that $\pi$ is a complete isometry, i.e., $F=0$.

\noindent
(d)${}\Rightarrow{}$(a)\enspace
Let $B$ be the hereditary \Csalg corresponding to $J$, i.e., $B=J^*\cap J$,
and $p_B$ the associated open projection in~$\Abar$, see Section~\ref{sect:inj-envel}.
Since, for each~$\alpha$, $0\leq e_\alpha\leq p_B$ we have $(1-p_B)e_\alpha=0$.
Hence $1-p_B\in\ell_{I(A)}(J)$ so, by assumption, $p_B=1$. As a result the
mapping $B\to I(A)$ is $\Oone$-essential which entails that $J\to I(A)$ has the
same property.

Assuming that $A$ is unital, the equivalence of (d) and~(e) is evident.
\end{proof}


The following proposition appears inter alia in~\cite{Arazth}.

\begin{proposition}\label{prop:algess}
Let\/ $J$ be an algebraically essential closed right ideal of the \C*~$A$.
Then\/ $J$ is\/ $\Oone$-essential.
\end{proposition}

For the proof, we need the following lemma \cite[Lemma~2.2]{Arazth}. We will make
use of the simple fact that $\ell(J)=0$ for every algebraically essential right
ideal~$J$. Indeed, if $x\in\ell(J)^*\cap J$ then $x^*x=0$; thus $\ell(J)^*\cap J=0$.
Since $\ell(J)^*$ is a right ideal, $J$ algebraically essential entails
that $\ell(J)=0$.

\begin{lemma}\label{lem:poscrit}
Let\/ $a\in A$ and\/ $J$ be an algebraically essential closed right ideal of the \C*~$A$.
Then\/ $a\geq0$ if and only if\/ $x^*ax\geq0$ for all\/ $x\in J$.
\end{lemma}

\begin{proof}
The ``only if''-part being evident we assume that $x^*ax\geq0$ for all $x\in J$.
Assume first that $a=a^*$. Then we can write $a=a_{+}-a_-$ with $a_+$ and $a_-$
positive and $a_+a_-=0$. Multiplying this identity by $a_-^{1/2}$ on the
right and on the left, we get $a_-^{1/2}aa_-^{1/2}=-a_-^2\le 0$. On
the other hand, $I=\{x\in A\mid a_-^{1/2}x\in J\}$ is an algebraically essential closed
right ideal of $A$, and for $x$ in $I$ we have, by assumption,
\begin{equation*}
0\le (a_-^{1/2}x)^*a(a_-^{1/2}x)=-x^*a_-^2x\le 0,
\end{equation*}
which entails $x^*a_-^2x=0$. Therefore $a_-\in\ell(I)=0$ implying that $a\geq0$
in this case.

An analogous argument shows that, if $b=b_+-b_-$ is self-adjoint and $x^*bx=0$
for all $x\in J$, then $b_+=b_-=0$ and so $b=0$. Let $a\in A$ be arbitrary with
$x^*ax\ge 0$ for all $x\in J$. Then $x^*(\frac{a-a^*}{2i})x=0$ for all $x\in J$ and
thus $a-a^*=0$ by the observation just made. Therefore $a=a^*$, and the first part of the
argument entails $a\ge 0$, as desired.
\end{proof}

As a consequence of this result we have
\begin{equation*}
\|a\|^2=\inf\{\lambda>0\mid x^*(\lambda^21-a^*a)x\geq0\ \,\forall\ x\in J\}
\qquad(a\in A).
\end{equation*}

\medskip\noindent
\begin{proofof}[Proposition~\ref{prop:algess}]\rm
Let $(e_\alpha)$ be an approximate identity of $J^*\cap J$. For $a\in A$,
take $\lambda\geq0$ such that $\|ae_\alpha\|\leq\lambda$ for all~$\alpha$.
Then, for all $\alpha$, $e_\alpha a^*ae_\alpha\leq\lambda^21$.
Therefore, for all $x\in J$, $(x^*e_\alpha)a^*a(e_\alpha x)\leq\lambda^2\,x^*x$
which implies that $x^*a^*ax\leq\lambda^2x^*x$.
By the above consequence of Lemma~\ref{lem:poscrit} it follows that
$\|a\|\leq\lambda$ wherefore
$\|a\|\leq\sup_\alpha\|ae_\alpha\|=\sup_\alpha\|e_\alpha a\|$.
Since the reverse inequality is obvious, we find $\|a\|=\sup_\alpha\|e_\alpha a\|$.

We can repeat this argument for matrices $(a_{kl})\in M_n(A)$ for any~$n$, since
$M_n(J)$ is an algebraically essential closed right ideal of $M_n(A)$ and
$\text{diag}(e_\alpha,\ldots,e_\alpha)$ is a left approximate identity
in $M_n(J)$. Hence, by applying criterion~(c) in Theorem~\ref{thm:ess-right-ideals}
we conclude that $J$ is $\Oone$-essential.
\end{proofof}

\medskip
By means of this we have established the first two of the implications between the
various notions of essentiality described above. The final one follows immediately
from the fact that, if $Jy=0$ for some $y\in A$ then $y^*e_\alpha=0$ for every
left approximate identity in the closed right ideal~$J$. Hence $y^*x=0$ for all $x\in J$,
wherefore $y=0$ provided $\ell_A(J)=0$.

\smallskip
Let $\BanA$ denote the category of Banach right $A$-modules with bounded
$A$-module maps as morphisms. The canonical choice for $\mathcal H$ in this case is the
class of all injective bounded right $A$-module maps. A closed right ideal $J$ of $A$ is
$\mathcal H$-essential (that is, $J\hookrightarrow A$ is $\mathcal H$-essential)
if and only if $J\cap K\ne0$ for every {\it closed\/} right ideal $K\ne0$ of~$A$.
The relation between essentiality of a closed right ideal $J$ in a \Calg $A$ in $\BanA$
and $\Oone$-essentiality is a priori unclear but for the fact that, if $J$ is
$\BanA$-essential, then $\ell_A(J)=0$. However, we have the following strong result.

For each $\epsilon >0$, we let $f_{\epsilon}$ denote the continuous real-valued function
on $[0,\infty)$ that is $0$ on the interval $[0,\epsilon]$, identically $1$ on the interval
$[\epsilon, 2\epsilon]$, and linear in the interval $[\epsilon,2\epsilon]$. Note that
$f_{\epsilon}(t)=th_{\epsilon}(t)$ for some continuous real-valued function $h_{\epsilon}$
and all $t\in[0,\infty)$.

\begin{lemma}\label{lem:baness-equal-algess}
Let\/ $J$ be a closed right ideal of a \Calg $A$. Then\/ $J$ is
algebraically essential if and only if\/ $J$ is\/ $\BanA$-essential.
\end{lemma}

\begin{proof}
Clearly, the ``only if''-part holds. Conversely, assume that $J$ is $\BanA$-essential
and let $a$ be a non-zero element in~$A$. Then there exists $\epsilon >0$ such that
$f_{\epsilon}(aa^*)\ne0$. Since $J$ is $\BanA$-essential, the closed right ideal
generated by $f_{\epsilon}(aa^*)$ has non-zero intersection with~$J$.
Let $z=\lim_{n\to\infty} f_{\epsilon}(aa^*)z_n$ be a non-zero element in~$J$.
Since $f_{\epsilon/2}(aa^*)f_{\epsilon}(aa^*)=f_{\epsilon}(aa^*)$, we see that
$z=f_{\epsilon /2}(aa^*)z$. On the other hand,  there is a positive element
$y=h_{\epsilon /2}(aa^*)$ in $A$ such that $f_{\epsilon/2}(aa^*)=aa^*y$, and
so $z=a(a^*yz)\in J\cap aA$ which entails that $J\cap aA\ne 0$. This proves that
$J$ is an algebraically essential right ideal of~$A$.
\end{proof}

Combining Proposition~\ref{prop:algess} and Lemma~\ref{lem:baness-equal-algess} we obtain
the following result.

\begin{proposition}\label{prop:baness-implies-o1ess}
Every\/ $\BanA$-essential closed right ideal of a \C*\/ $A$ is $\Oone$-essential.
\end{proposition}

The converse of the statement in Proposition~\ref{prop:baness-implies-o1ess}
does not hold; in Section~\ref{sect:example} we shall provide an example
(Example~\ref{exam:o1ess-not-alg-ess}).

We will make frequent use of the following result.

\begin{lemma}\label{lem:cbnorm-equal-norm}
Let\/ $J$ be a closed right ideal of a \Calg $A$, and let\/ $g\colon J\to M$
be an\/ $A$-module map into an operator right\/ $A$-module~$M$. If\/ $g$ is
bounded then\/ $g$ is completely bounded and\/
$\cbn g=\|g\|=\sup_\alpha\|g(e_\alpha)\|$ for every approximate identity\/
$(e_\alpha)$ in\/ $J^*\cap J$.
\end{lemma}

\begin{proof}
Let $x\in J$, $\|x\|=1$, and suppose that $g$ is bounded. Then
\[
\|g(x)\|=\|\lim_\alpha g(e_\alpha x)\|\leq\sup_\alpha\|g(e_\alpha)\|,
\]
whence $\|g\|=\sup_\alpha\|g(e_\alpha)\|$. Noting that $\text{diag}(e_\alpha,\ldots,e_\alpha)$
is a left approximate identity for $M_n(J)$, we can use the same argument to find
$\|g^{(n)}\|=\sup_\alpha\|g(e_\alpha)\|=\|g\|$ for every $n\in\NN$.
\end{proof}

We are next going to study the relationship between open projections in the regular monotone
completion $\Abar$ of a \Calg $A$ and certain one-sided ideals of~$A$.

A closed right ideal $J$ of $A$ will be called {\em $\Oone$-essentially closed\/} in case
$J$ has no proper $\Oone$-essential extension in~$A$; cf.~Theorem~\ref{thm:ess-right-ideals}.
Denote by $\mathfrak G(A)$ the set of those right ideals.

\begin{lemma}\label{lem:openproj}
There is a bijective correspondence between\/ $\mathfrak G(A)$ and the set of
open projections in\/ $\Abar$ given by the rules
$$J\mapsto{\Abarsa}\hbox{-}\sup e_\alpha,\qquad  p\mapsto p\Abar\cap A,$$
where\/ $(e_\alpha)$ is any approximate identity in\/ $J\cap J^*$.
\end{lemma}

\begin{proof}
Let $J$ be an $\Oone$-essentially closed right ideal of $A$ and set
$p={\Abarsa}\hbox{-}\sup e_\alpha$, where\/ $(e_\alpha)$ is any approximate identity
in $J\cap J^*$.
Using Theorem~\ref{thm:ess-right-ideals} it is easily seen that $pI(A)$ is the
injective envelope of $J$ in the category $\COoneA$ and that $J$ is essential in
$pI(A)\cap A$ in this category. Since $J$ has no proper $\Oone$-essential extension in~$A$,
it follows that $pI(A)\cap A=p\Abar\cap A=J$.

Conversely, let $p$ be an open projection in $\Abar$, so that
$p={\Abarsa}\hbox{-}\sup a_{\beta}$ for an increasing net
$(a_{\beta})$ of positive elements in~$A$. Set $J=p\Abar \cap A$.
If $(e_\alpha)$ is an approximate identity in $J\cap J^*$, then
necessarily we have $p={\Abarsa}\hbox{-}\sup e_\alpha$. Indeed,
suppose that ${\Abarsa}\hbox{-}\sup e_\alpha=q\le p$. Since
$a_{\beta}=\lim_\alpha e_\alpha a_{\beta}$, we get
$qa_{\beta}=a_{\beta}$, so that $p={\Abarsa}\hbox{-}\sup a_{\beta}\le q$.
\end{proof}

Thus, for every closed right ideal $J$, we can find a
maximal $\Oone$-essential extension of $J$ in $A$, namely
$I(J)\cap A=p\Abar\cap A$, where $p={\Abarsa}\hbox{-}\sup e_\alpha$ for
an approximate identity $(e_\alpha)$ of $J\cap J^*$.

We introduce an important class of closed one-sided ideals, the
annihilator ideals. Denote by $\mathfrak A_r(A)$ (resp.,
$\mathfrak A_\ell(A)$) the set of all right (resp., left) annihilators
in $A$ of subsets in~$A$; cf.~page~\pageref{p:annis}. Observe that $\mathfrak A_r(A)$ is
indeed the set of right annihilators of (closed) left ideals of~$A$
and there exists a bijective correspondence $\mathfrak A_r(A)\to
\mathfrak A_\ell(A)$ given by taking left annihilators, which intertwines the involution:
$\ell(J^*)=r(J)^*$ for each $J\in\mathfrak A_\ell(A)$.

Let $J$ be a closed right ideal of $A$ and let $J'=p\Abar \cap A$
be its maximal $\Oone$-essential extension in $A$. Then
$\ell_{\smallAbar}(J)=\ell_{\smallAbar} (J')=\Abar(1-p)$, so in particular we
have $\ell _A(J)=\ell _A(J')=\Abar (1-p)\cap A$; thus
\[
J\subseteq J'\subseteq r_A\ell_A(J).
\]
If $J\in\mathfrak A_r(A)$ then $J=r_A\ell _A(J)$ and, in particular,
$J=J'$; therefore $J$ does not admit any proper $\Oone$-essential extension in~$A$.
This shows that $\mathfrak A_r(A)\subseteq \mathfrak G (A)$,
and the set of right annihilator ideals can be identified with a
subset of open projections in~$\Abar$.

Another family of closed right ideals that we can consider is the
family of right ideals of the form $p\Abar \cap A$, where $p$ is a
clopen projection in $\Abar$, that is, both $p$ and $1-p$ are open
projections. Let us denote the set of such right ideals by $\mathfrak F(A)$.\label{p:FA}

\begin{proposition}\label{prop:clopenproj}
For a \C*\/ $A$ we have the following inclusions
$$\mathfrak F(A)\subseteq \mathfrak A_r(A)\subseteq \mathfrak G(A).$$
Moreover, a right ideal\/ $J$ is in\/ $\mathfrak F(A)$ if and only if
it does not admit any proper\/ $\Oone$-essential extension and\/ $J\oplus r_A(J^*)$
is an\/ $\Oone$-essential closed right ideal of~$A$.
\end{proposition}

\begin{proof}
We have already observed that $\mathfrak A_r(A)\subseteq \mathfrak G(A)$.
Let $p$ be a clopen projection in~$\Abar$. Then $1-p$ is
the supremum in $\Abar$ of an increasing net of positive elements
of $A$, and so $r _A(\Abar (1-p)\cap A)=p\Abar \cap A$, which
shows that $p\Abar \cap A\in \mathfrak A_r(A)$.

Let $J=p\Abar \cap A$ be an $\Oone$-essentially closed right ideal of $A$,
where $p$ is an open projection in $\Abar$. Let $(e_\alpha)$ be
an approximate identity in $J\cap J^ *$ and let $(f_\beta)$ be an
approximate identity for $\ell_A(J)\cap\ell _A(J)^*$. Then the
supremum in $\Abar $ of $(f_\beta)$ is $1-p$ if and only if the
supremum in $\Abar$ of $(e_\alpha+f_\beta)$ is $1$ which in turn is equivalent to
\hbox{$J\oplus r_A(J^*)$} being a closed $\Oone$-essential right ideal of~$A$. This
shows the last statement.
\end{proof}

\section{The Maximal {\sl C*}-Algebra of Quotients}\label{sect:max-quotients}

In this section we shall first of all recall the concept of the maximal \C* of quotients.
It is the analytic analogue of the maximal symmetric ring of quotients studied by
Lanning~\cite{Lann}, much in the same way as the local multiplier algebra, see
\cite{AMMzth}, is the analytic counterpart of the symmetric ring of quotients.

Following the approach used in \cite{Arazth}, we let $A$ be a unital \Calg and first
construct the maximal symmetric algebra of quotients, $\qmaxs$, of~$A$.
Let $\Ier$ denote the filter of algebraically essential right ideals of~$A$.
We consider triples $(f,g,I)$, where $I\in\Ier$, $f\colon I^*\rightarrow A$ is a left
$A$-module homomorphism and $g\colon I\rightarrow A$ is a right $A$-module homomorphism
satisfying the compatibility rule
\[
f(x)y=xg(y)\qquad(x\in I^*,\,y\in I).
\]
Two such triples $(f_1,g_1,I_1)$ and $(f_2,g_2,I_2)$ are said to be {\it equivalent\/}
if $g_1$ and $g_2$ coincide on $I_1\cap I_2$. It immediately follows that $f_1$ and
$f_2$ agree on $I_1^*\cap I_2^*$, which also shows that the first component in the
triple $(f,g,I)$
is uniquely determined by the second; thus the {\it existence\/} of~$f$ is the crucial
assumption. Let $\qmaxs$ be the set of all equivalence classes of triples $(f,g,I)$
of this kind. We define algebraic operations on $\qmaxs$ by
\begin{equation*}
\begin{split}
[(f_1,g_1,I_1)]+[(f_2,g_2,I_2)] &=[(h_1,h_2,I_1\cap I_2)],\\
[(f_1,g_1,I_1)][(f_2,g_2,I_2)]  &=[(k_1,k_2,J)],\\
[(f,g,I)]^*                     &=[(g^*,f^*,I)].
\end{split}
\end{equation*}
The right hand sides are defined as follows:
$h_1(x)=f_1(x)+f_2(x)$ and $h_2(y)=g_1(y)+g_2(y)$;
$J=f_1^{-1}(I_2^*)^*\cap g_2^{-1}(I_1)$ and $k_1(x)=f_2(f_1(x))$,
$k_2(y)=g_1(g_2(y))$;
$f^*(x)=f(x^*)^*$, $x\in I$ and $g^*(x)=g(x^*)^*$, $x\in I^*$. Endowed with these
operations, $\qmaxs$ becomes a unital complex *-algebra called the
{\it maximal symmetric algebra of quotients of\/}~$A$, cf.~\cite{Lann}.
The \C* $A$ is canonically embedded into $\qmaxs$ via
$a\mapsto(R_a,L_a,A)$, where $R_a$ (resp.\ $L_a$) denotes right (resp.\ left)
multiplication by~$a\in A$.

In order to define a `bounded part' of $\qmaxs$, and a corresponding {\sl C*}-norm
on it, we first introduce a positive cone, $\qmaxs_+$, of~$\qmaxs$. We say that
$c\in\qmaxs$ is {\it positive\/} if, for some representative $(f,g,I)$ of~$c$,
\[
f(x^*)x=x^*g(x)\geq0\qquad(x\in I).
\]
Lemma~\ref{lem:poscrit} implies that $A_+=A\cap\qmaxs_+$.
It is easily verified that $\qmaxs_+$ is indeed a convex cone \cite[p.~15]{Arazth}.

An element $d\in\qmaxs$ is said to be {\it bounded\/} if, for some $\lambda\geq0$,
both $\lambda^2\,1-dd^*$ and $\lambda^2\,1-d^*d$ belong to~$\qmaxs_+$.
We shall denote the set of all bounded elements in $\qmaxs$ by~$\qmaxsb$ and we
introduce a norm by setting
\[
\|d\|=\inf\bigl\{\lambda\geq0\mid \lambda^2\,1-dd^*,\ \lambda^2\,1-d^*d\in\qmaxs_+\bigr\},
\]
where $d\in\qmaxsb$.

The following result is proved in \cite[Theorem~2.3]{Arazth}. Note that, if $A$
is non-unital, we apply the above construction to its minimal unitisation.

\begin{theorem}\label{thm:bdd-part-qmaxs}
For every \C*\/ $A$, $\qmaxsb$ endowed with the above norm is a unital pre-\C*
containing\/ $A$ isometrically.
\end{theorem}

\begin{definition}
The completion of $\bigl(\qmaxsb,\|\cdot\|\bigr)$ is called {\it the maximal
\C* of quotients of\/} $A$ and will be denoted by~$\qmax$.
\end{definition}

\begin{remark}\label{rem:ess-right-ideals-vs-algebras}
Every element $d\in\qmaxsb$ has representatives of the form $(f,g,J)$,
where both $f$ and $g$ are bounded maps on a closed essential right ideal~$J$,
and conversely, $d\in\qmaxs$ is bounded if it has a representative of this form.
Let $\Icer$ denote the filter of essential closed right ideals in the \Calg~$A$.
(Note that, by Lemma~\ref{lem:baness-equal-algess} above, for a closed right ideal $J$
the conditions to be algebraically essential and essential in $\BanA$ are equivalent
and we shall thus unambiguously call such an ideal ``essential".)
We make use of the well-known bijective correspondence between closed right ideals $J$
and hereditary \Cs*s $D$ of $A$ given by
$D=J^*\cap J$ and $J=\{x\in A\mid xx^*\in D\}$, see, e.g., \cite[Theorem~1.5.2]{Pedsvn}.
We call the hereditary \Cs* $D$ {\it essential\/} if its corresponding closed right
ideal $J$ is essential and shall denote the set of all essential hereditary \Csalgs
by~$\Heress$. If $d\in\qmaxsb$ and $(f,g,J)$ is a representative as above, then
$\|d\|=\|f\|=\|g\|$ by \cite[Lemma~2.4]{Arazth}. In particular, every multiplier
$c\in M(D)$, where $D\in\Heress$, gives rise to an element $[(R_c,L_c,J)]$ in
$\qmaxsb$ (where $J$ is the closed right ideal generated by~$D$) with the same norm as~$c$.
\end{remark}

In order to understand the relation between the multiplier algebras of essential
hereditary \Csalgs and $\qmax$ indicated in
Remark~\ref{rem:ess-right-ideals-vs-algebras} above more closely, we need the following
universal property of the maximal \Calg of quotients.

Whenever $B$ is a \C* containing a *-isomorphic image of $A$, we call $B$
{\it an enlargement\/} of~$A$. In this case, we write $A\hookrightarrow B$ and identify
$A$ with its image in~$B$. For each $J\in\Icer$, we define an $A$-subbimodule $B_J^A$ of~$B$
by
\[
B_J^A=\{b\in B\mid bJ+b^*J\subseteq A\}.
\]
Set $B_\Icer=\bigcup_{J\in\Icer} B_J^A$. Then $B_\Icer$ is a *-subalgebra of $B$
so its completion is a \Calgc. We say that $B$ is {\it an $\Icer$-enlargement\/}
if $B_\Icer$ is dense in~$B$. Note that $A\hookrightarrow\qmax$ is an
$\Icer$-enlargement.

\begin{theorem}\label{thm:univ-prop}{\rm(Universal Property of $\qmax$)}\quad
Let\/ $A\hookrightarrow B$ be an enlargement of the \C*~$A$.
Then there exists a unique contractive *-homomorphism\/ $\psi\colon B_\Icer\to\qmax$
which is the identity on~$A$. The mapping\/ $\psi$ is isometric if and only if\/
$B$ is an $\Icer$-enlargement such that, whenever\/ $b\in B$ and\/ $J\in\Icer$,
$bJ=0$ implies that\/ $b=0$.
\end{theorem}

This property, which is completely analogous to the universal property of the
local multiplier algebra, see \cite[Proposition~2.3.4]{AMMzth}, is obtained in
\cite[Proposition~2.5]{Arazth}.

The local multiplier algebra $\mloc$ of a \Calg $A$ is defined by
$\mloc=\Dirlim{\Ice}\,M(I)$, where $I$ runs through the filter $\Ice$ of closed essential
(two-sided) ideals in~$A$; see \cite[Section~2]{AMMzth}. As a direct consequence of
this definition, the local multiplier algebra embeds canonically into the maximal
\Calg of quotients.

\begin{proposition}\label{prop:mloc-into-qmax}
For every \Calg $A$, there is a unique *-isomorphism from\/ $\mloc$ into\/
$\qmax$ which is the identity on~$A$.
\end{proposition}

\begin{proof}
As observed in Remark~\ref{rem:ess-right-ideals-vs-algebras} above, for every
closed essential ideal $I$ in $A$, the multiplier algebra $M(I)$ embeds isometrically
into~$\qmax$. If $J\in\Ice$, $J\subseteq I$ then $M(I)$ embeds canonically into $M(J)$
by restriction of the multipliers (these are the connecting maps in the direct limit
construction). Since these embeddings are clearly compatible with each other,
we obtain a *-isomorphism from $\mloc$ into $\qmax$ which is the identity on~$A$.
The uniqueness is evident.
\end{proof}

It is also clear that the embedding $\mloc\to\qmax$ is precisely the mapping $\psi$
in Theorem~\ref{thm:univ-prop} and that $\mloc$ is an $\Icer$-enlargement of the
above kind. The situation is more complicated once we look at arbitrary essential
hereditary \Cs*s.

\begin{example}\label{exam:hered-subalgs}
For $D\in\Heress$, we consider $M(D)$ as a \Cs* of $\qmax$ as described in
Remark~\ref{rem:ess-right-ideals-vs-algebras} above. Let $B$ be the \C* generated
by $\bigcup_{D\in\Heress}M(D)$. If $J$ is the essential closed right ideal corresponding
to $D\in\Heress$, then $M(D)\subseteq B_J^A$; hence the *-subalgebra generated by
$\bigcup_{D\in\Heress}M(D)$ is contained in $B_\Icer$. It follows that $B_\Icer$ is
dense in $B$ wherefore $A\hookrightarrow B$ is an $\Icer$-enlargement of $A$ and
the mapping $\psi$ described in Theorem~\ref{thm:univ-prop} is simply the inclusion
$B\subseteq\qmax$.

A description of the \C* $B$ as a direct limit of the multiplier algebras
$M(D)$, $D\in\Heress$ is not available.
\end{example}

We now turn our attention to the relation between $\qmax$ and the injective envelope
$I(A)$ of a \C*~$A$.

The following result can be extracted from \cite{FrankPaulsen} but we include a proof
for the sake of completeness.

\begin{lemma}\label{lem:multipliers-in-injenv}
Let\/ $A$ be a unital \C*, and let\/ $J$ be an essential closed right ideal of~$A$.
For every bounded right $A$-module homomorphism\/ $g\colon J\to A$ there is a unique
element\/ $y\in I(A)$ such that\/ $g=L_y$ and\/ $\|g\|=\|y\|$.
\end{lemma}

\begin{proof}
By Lemma~\ref{lem:cbnorm-equal-norm}, $g$ is completely bounded and thus we can assume,
without restricting the generality, that $\cbn g=1$. By the module version of
Wittstock's extension theorem, $I(A)$ is injective in $\COoneA$; hence there exists
a completely contractive right $A$-module map $\tilde g\colon A\to I(A)$ such
that $\rest{\tilde g}{J}=g$ and $\cbn{\tilde g}=1$. From
$\tilde g(a)=\tilde g(1\!\cdot\!a)=\tilde g(1)a$
for all $a\in A$ we see that $\tilde g$ is left multiplication by some element
$y\in I(A)$. Suppose that $y'\in I(A)$ is another element with the property
$g(x)=y'x$ for all $x\in J$. Then $y-y'\in\ell_{I(A)}(J)$ so that $y-y'=0$
by Theorem~\ref{thm:ess-right-ideals}. As a result there is a unique element $y\in I(A)$
such that $g=L_y$ and $\|g\|=\cbn g=\|y\|$.
\end{proof}

Invoking the universal property of $\qmax$ (Theorem~\ref{thm:univ-prop})
we have a unique contractive *-homomorphism
\[
\psi\colon I(A)_\Icer\to\qmax
\]
which is the identity on~$A$. Let $d\in\qmaxsb$ and $(f,g,J)$ be a representative
of $d$ such that $J\in\Icer$. By Lemma~\ref{lem:multipliers-in-injenv},
there is a unique element $y\in I(A)$ such that $g=L_y$ and $\|g\|=\|y\|$.
Since $(f(x)-xy)z=x(g(z)-yz)=0$ for all $x\in J^*$ and $z\in J$,
and $\ell_{I(A)}(J)=0$ by Proposition~\ref{prop:baness-implies-o1ess}
and Theorem~\ref{thm:ess-right-ideals}, it follows that $f=R_y$.
If $(f',g',J')$ is another such representative of~$d$, then
$g'=g=L_y$ on $J'\cap J$ and hence $f'=f=R_y$ on $J'\cap J$. Consequently,
\[
\tilde\psi\colon\qmaxsb\to I(A)_\Icer,\quad d\mapsto y
\]
is a well-defined contractive *-homomorphism which, by construction, is
a left inverse of~$\psi$.

We therefore obtain the following result, which is the first part 
of \cite[Theorem~2.7]{Arazth}.

\begin{theorem}\label{thm:qmax-into-injenv}\quad
For every \C* $A$, the maximal \Calg of quotients\/ $\qmax$ embeds
canonically into the injective envelope\/ $I(A)$ of $A$ such that\/
$\qmaxsb$ is isometrically *-isomorphic to\/ $I(A)_\Icer$.
\end{theorem}

Indeed, it is fairly easy to show that $\qmax$ is contained in the regular monotone
completion $\Abar$ of the \C*~$A$.

\begin{proposition}\label{prop:qmax_into_abar}
For every \C*\/ $A$, we have\/ $\qmax\subseteq\Abar$.
\end{proposition}

\begin{proof}
It suffices to show that every positive element $y$ in $\qmaxsb\subseteq I(A)$
belongs to~${\Abar}$. By Theorem~\ref{thm:qmax-into-injenv}, there is a closed essential
right ideal $I$ in $A$ such that $yI\subseteq A$. Let $(u_\lambda)$ be an
approximate identity in~$I^*\cap I$. Then ${I(A)_{sa}}$-$\sup_\lambda u_\lambda=1$;
\cite[Lemma~1.9 and Theorem~3.1]{Hameoa} therefore entail that
\[
y^2=y\cdot{I(A)_{sa}}\text{-}\sup_\lambda u_\lambda\cdot y
   ={I(A)_{sa}}\text{-}\sup_\lambda\,yu_\lambda y
   ={\Abar_{\!sa}}\text{-}\sup_\lambda\,yu_\lambda y\in\Abar.
\]
As a result $y$ belongs to~$\smash{\Abar}$.
\end{proof}

\begin{corollary}\label{cor:qmax_for_mon_complete}
For a monotone complete \C*\/ $A$, we have\/ $\qmax=A$.
\end{corollary}

\begin{corollary}\label{cor:qmax_vneumann}
The maximal \C* of quotients of a von Neumann algebra coincides with the algebra.
\end{corollary}

In the next result we give a detailed proof of the second part of \cite[Theorem~2.7]{Arazth}.

\begin{theorem}\label{theor:centres}
Let\/ $A$ be a \C*. Then
\[
Z(\mloc)=Z(\qmax)=Z(\Abar)=Z(I(A)).
\]
\end{theorem}

\begin{proof}
By \cite[Theorem 6.3]{Hameoa} and Proposition~\ref{prop:invariant-injective},
we have $Z(\Abar)=Z(I(A))$. Hence, by \cite[Corollary~4.3]{Hamsvn}, we conclude that
\begin{equation*}\label{equ:Open1}
Z(\Mloc{A})\subseteq Z(\Qmax{A})\subseteq Z(\Abar)=Z(I(A)).
\end{equation*}
Let $p$ be a central projection in~$\Abar$. It follows from
\cite[Lemma~1.3]{Hameta} that $p$ is clopen and that the closed ideal
$I=p\Abar\cap A$ satisfies $I=I^{\perp\perp}$. Moreover,
$p$ is the obvious projection in $\Mloc{A}$ given by projecting
the essential ideal $I\oplus I^{\perp}$ onto~$I$. Consequently, all the
projections in $Z(\Abar)$ belong to $Z(\Mloc{A})$, and since
$Z(\Abar)$ is an \AW*, we obtain that $Z(\Mloc{A})=Z(\Abar)$. Thus the above inclusions
are in fact equalities.
\end{proof}

\begin{remark}\label{rem:inclus-sets-of-ideals}
Let $\mathfrak A(A)$ be the set of closed ideals $I$ in $A$ satisfying $I=I^{\perp\perp}$.
By the arguments in the proof of Theorem~\ref{theor:centres}, we have a bijective
correspondence between $\mathfrak A(A)$ and the set of projections in $Z(I(A))$
given by $I=p_I\Abar\cap A$. The projections in $\qmaxsb$ correspond to $\Oone$-essentially
closed right ideals $J$ of $A$ such that $J\oplus r_A(J^*)$ is algebraically essential
in~$A$. By Proposition~\ref{prop:algess}, the set $\mathfrak F_a(A)$ of these right ideals is
contained in~$\mathfrak F(A)$, the set of right ideals introduced on
page~\pageref{p:FA}. Hence, we obtain the following chain of inclusions
\[
\mathfrak A(A)\subseteq\mathfrak F_a(A)\subseteq\,\mathfrak F(A)\subseteq
\mathfrak A_r(A)\subseteq\mathfrak G(A).
\]
Each of these sets corresponds to a certain subset of open projections in~$\Abar$.
\end{remark}

We proceed to use the results in Section~\ref{sect:ess-right-ideals}\ in order
to discuss the operator algebra structure of the maximal \C* of quotients.

For a \C* $A$, let $\COA$ denote the category whose objects are the operator
right $A$-modules and whose morphisms are the completely bounded right $A$-module
maps. We simply denote the space of all morphisms between two operator right $A$-modules
$E$ and $F$ by $\CBsA(E,F)$. Note that, if $F$ is an operator $A$-bimodule,
this space is an operator left $A$-module, where
the operator space structure is given by $M_n(\CBsA(E,F))=\CBsA(E,M_n(F))$, $n\in\NN$
and the $A$-module structure on $\CBsA(E,F)$ is defined by
$(ag)(x)=ag(x)$, for $a\in A$, $x\in E$ and $g\in\CBsA(E,F)$.

Let $A$ be a unital \C*, and let $I$, $J$ belong to $\Icer$. Denote by
$\rho_{JI}$ the restriction maps
$\rho_{JI}\colon\CBsA(I,A)\to\CBsA(J,A)$, whenever $J\subseteq I$.
Let $E_b(A)=\Alglim{\Icer}\,\CBsA(I,A)$; we will now see that this is in fact
a direct limit in the category of operator left $A$-modules.

\begin{proposition}\label{prop:mapstau}
Let\/ $A$ be a unital \C*. For every\/ $I\in\Icer$ there exists a completely
isometric left $A$-module homomorphism\/ $\tau_I\colon\CBsA(I,A)\to I(A)$.
Moreover, if\/ $I,\,J\in\Icer$ and\/
$J\subseteq I$ then\/ $\tau_I=\tau_J\circ \rho_{JI}$, where\/
$\rho_{JI}\colon\CBsA(I,A)\to\CBsA(J,A)$ is completely isometric.
\end{proposition}

\begin{proof}
Let $I$ be a closed essential right ideal of~$A$, and let $g\colon I\to A$
be a (completely) bounded right $A$-module map (cf.\ Lemma~\ref{lem:cbnorm-equal-norm}).
By Lemma~\ref{lem:multipliers-in-injenv}, there exists a
unique element $y\in I(A)$ such that $g=L_y$ and
$\|y\|=\|g\|=\cbn g=\sup_\alpha\|g(e_\alpha)\|$, where $(e_\alpha)$ is some
left approximate identity in~$I$. By means of this, we
obtain a well-defined $A$-linear isometry $\tau_I\colon\CBsA(I,A)\to I(A)$ which
evidently satisfies the compatibility condition $\tau_I=\tau_J\circ \rho_{JI}$.

In order to show that $\tau_I$ is completely isometric, let $n\in\NN$ and
consider
\[
\tau_I^{(n)}\colon\CBsA(I,M_n(A))\to I(M_n(A))=M_n(I(A)).
\]
Take $\hat g=(g_{ij})\in\CBsA(I,M_n(A))$; then
$\tau_I^{(n)}(\hat g)=(\tau_I(g_{ij}))=(y_{ij})$,
where $y_{ij}\in I(A)$ is the unique element such that $g_{ij}=L_{y_{ij}}$.
By Lemma~\ref{lem:cbnorm-equal-norm} we have
\begin{equation*}
\|\hat g\|=\sup_\alpha\|\hat g(e_\alpha)\|
          =\sup_\alpha\|(g_{ij}(e_\alpha))\|
          =\sup_\alpha\|(y_{ij}e_\alpha)\|
          \leq\|(y_{ij})\|.
\end{equation*}
On the other hand, letting $\hat y=(y_{ij})\in M_n(I(A))$,
$L_{\hat y}\colon M_n(I)\to M_n(A)$ defines a bounded right $M_n(A)$-module
homomorphism on the closed essential right ideal $M_n(I)$ of $M_n(A)$.
Applying Lemma~\ref{lem:cbnorm-equal-norm} with the approximate identity
$\text{diag}(e_\alpha,\ldots,e_\alpha)$ we find that
$\|L_{\hat y}\|=\sup_\alpha\|(y_{ij}e_\alpha)\|$ and since $M_n(I(A))$ is injective,
it follows that $\|\hat y\|=\|L_{\hat y}\|$ (Lemma~\ref{lem:multipliers-in-injenv}).
As a result $\|\hat y\|=\|\tau_I^{(n)}(\hat g)\|=\|\hat g\|$, wherefore
$\tau_I$ is a complete isometry.

Finally,
\[
\bigl\|\rho_{JI}^{(n)}(\hat g)\bigr\|
=\bigl\|\tau_J^{(n)}\bigl(\rho_{JI}^{(n)}(\hat g)\bigr)\bigr\|
=\bigl\|\tau_I^{(n)}(\hat g)\bigr\|=\|\hat g\|
\]
shows that each $\rho_{JI}$ is completely isometric.
\end{proof}

Hence the algebraic direct limit $E_b(A)$ has a canonical operator space structure.
Indeed, $E_b(A)$ is a unital operator algebra with the product given by the formula
\[
[(f,I)]\,[(g,J)]=[(fg,g^{-1}(I))],
\]
where $[(f,I)]$ denotes the equivalence class of $f\in\CBsA(I,A)$, $I\in\Icer$
in~$E_b(A)$.

\begin{corollary}\label{cor:eb-into-injenv}
Let\/ $A$ be a unital \Calgc, and let\/ $E_b(A)$ be the operator algebra constructed above.
Then there is a completely isometric multiplicative embedding\/ $E_b\to I(A)$
which is the identity on~$A$.
\end{corollary}

\begin{proof}
This follows from Proposition~\ref{prop:mapstau} and the above remarks.
\end{proof}

Now we can give a new description of the maximal \Calg of quotients.
Considering $E_b(A)$ as an operator subalgebra in $I(A)$, we have that
$\qmaxsb=E_b(A)\cap E_b(A)^*$ and its completion in~$I(A)$ is~$\qmax$.

\section{The Case of {\sl AW*}-Algebras}\label{sect:awstars}

An outstanding question is whether every \AW* is monotone complete;
if the answer is positive, Corollary~\ref{cor:qmax_for_mon_complete} entails
that the maximal \C* of quotients of every \AW* agrees with the algebra itself.
In this section we shall employ different methods to obtain the latter statement
for some subclasses of \AW*s.

We first consider \AW*s with the property that each family of mutually orthogonal,
non-zero projections is countable. Such algebras will be called {\it\sigfinite/\/},
and we shall see that they have some interesting properties.
First, we will review a couple of known results.

\begin{lemma}\label{lem:principal-ann}
Let\/ $A$ be a \sigfinite/ \AW*. Then every closed right ideal\/ $I$ of\/ $A$
contains a principal closed right ideal\/ $J$ such that $\ell(I)=\ell(J)$.
\end{lemma}

\begin{proof}
By \cite[Proposition~1.5]{Wrighteza}, there exists a countable family $\{e_n\mid n\in\NN\}$
of projections in $I$ such that $\bigvee_n e_n=\bigvee P(I)$, where $P(I)$ denotes
the set of all projections in~$I$. Let $a=\sum_{n=1}^\infty 2^{-n}e_n$ and put
$J=\ol{aA}$. Since $\ell(\{a\})=\ell(\{e_n\mid n\in\NN\})$ it follows that
$\ell(J)=\ell(I)$.
\end{proof}

We include the proof of the following well-known result for completeness.

\begin{lemma}\label{lem:principal-hom}
Let\/ $A$ be a \C*, and let\/ $a,b\in A$. The mapping\/ $f(a)=b$ extends
to a bounded $A$-module homomorphism from\/ $\ol{aA}$ into\/ $A$ if and only if
there exists\/ $\gamma>0$ such that\/ $b^*b\leq \gamma\,a^*a$.
\end{lemma}

\begin{proof}
Suppose first that $b^*b\leq \gamma\,a^*a$ for some $\gamma>0$. Then, for each $x\in A$,
$\|bx\|^2\leq\gamma\,\|ax\|^2$ implying that $f\colon aA\to A$, $ax\mapsto bx$ is a
well-defined bounded $A$-module homomorphism with $\|f\|\leq\gamma^{1/2}$
and thus extends to~$\ol{aA}$.

Conversely, suppose that $f\colon\ol{aA}\to A$ given by $f(ax)=bx$ for all $x\in A$
is bounded.
Let $\eps>0$ and put $y=(\gamma+\eps)\,a^*a-b^*b$, where $\gamma=\|f\|^2$.
Then, for $x=y^{-}$, we have
\[
(\gamma+\eps)\,xa^*ax-xb^*bx=xyx=-x^3.
\]
Therefore $\|bx\|^2\geq(\gamma+\eps)\,\|ax\|^2>\gamma\,\|ax\|^2$ unless $ax=0$.
Since $\|bx\|=\|f(ax)\|\leq\|f\|\,\|ax\|$, it follows that we must have $ax=0$.
Hence $bx=0$ and thus $x=0$. This shows that $b^*b\leq(\gamma+\eps)\,a^*a$ for every
$\eps>0$ which entails that $b^*b\leq\gamma\,a^*a$, since $A_+$ is closed.
\end{proof}

\begin{proposition}\label{prop:qmaxs-sigf-aw}
Let\/ $A$ be a \sigfinite/ \AW*. Then\/ $\qmaxsb=A$. More precisely, for every closed right
ideal\/ $I$ of\/ $A$ with\/ $\ell_A(I)=0$ and\/ each bounded right $A$-module map\/
$f\colon I\to A$, there is a unique element\/ $c\in A$ such that\/ $f=L_c$.
\end{proposition}

\begin{proof}
Let $f\colon I\to A$ be a bounded right $A$-module map.
We shall first show that there is a unique element $c\in A$ such that
$f(e)=ce$ for all projections $e$ in~$I$.

Let $\{e_n\mid n\in\NN\}$ be a countable family of projections in $I$ as in
Lemma~\ref{lem:principal-ann}, that is, with the property
$\ell_A\bigl(\{e_n\mid n\in\NN\}\bigr)=\ell_A(I)=0$.
Set $a=\sum_{n=1}^\infty 2^{-n}e_n$; then $e_n\leq2^n\, a$ for all~$n$ and $\ell_A(\{a\})=0$.
Note also that $\ol{aA}=\ol{\sum_n e_nA}$.

Let $b=f(a)$; by Lemma~\ref{lem:principal-hom} we have $b^*b\leq \gamma\,a^*a$ for some
$\gamma>0$. Since every \AW* is a {\sl UMF}-algebra by \cite[Corollary~3.6]{AGoldnth},
there is $c\in A$ such that $b=ca$ and $\|c\|\leq \gamma^{1/2}$. It follows that
$\rest{f}{\ol{aA}}=L_c$, in particular $f(e_n)=ce_n$ for all $n\in\NN$.

Take any two projections $e$ and $e'$ in~$I$. By applying the above procedure to
both a countable family $P$ of projections in $I$ containing~$e$ and a countable
family $P'$ of projections in $I$ containing~$e'$ we obtain $c,\,c'\in A$
such that $f(p)=cp$ for all $p\in P$ and $f(p')=c'p'$ for all $p'\in P'$,
respectively. By applying this procedure to $P\cup P'$ we obtain yet another element
$c''\in A$ such that $f(p'')=c''p''$ for all $p''\in P\cup P'$. It follows that
\[
(c-c'')p=0\text{ for all }p\in P\quad\text{and}\quad
(c'-c'')p'=0\text{ for all }p'\in P'.
\]
Since $\ell_A(P)=0=\ell_A(P')$, we conclude that $c=c''=c'$. Hence
$f(e)=ce$ for all projections $e$ in~$I$ and $c$ is the only element in $A$ with
this property.

For an arbitrary element $x\in I$ and $\eps>0$, let $e\in I$ be a projection
such that $\|x-ex\|<\eps$ \cite[Corollary on~p.~43]{Berb}.
Since by the above $f(ex)=f(e)x=cex$, we conclude that
\begin{equation*}
\|f(x)-cx\|\leq\|f\|\,\|x-ex\|+\|c\|\,\|ex-x\|<(\|f\|+\|c\|)\,\eps,
\end{equation*}
wherefore $f(x)=cx$ as claimed.
\end{proof}

\begin{corollary}\label{cor:qmax-sigf-aw}
For every \sigfinite/ \AW*\/ $A$ we have\/ $\qmax=A$.
\end{corollary}

Call a \C* $A$ {\it$A$-injective\/} if every bounded right $A$-module homomorphism
from a closed right ideal in $A$ into $A$ is given by left multiplication by an element
in~$A$. Proposition~\ref{prop:qmaxs-sigf-aw} enables us to obtain the following result.

\begin{proposition}\label{prop:sigf-aw-selfinj}
Every \sigfinite/ \AW*\/ $A$ is\/ $A$-injective.
\end{proposition}

\begin{proof}
Let $f\colon I\to A$ be a bounded right $A$-module homomorphism defined on the closed
right ideal $I$ of~$A$. Put $J=I+\ell_A(I)^*=I+r_A(I^*)$; this is a right ideal of $A$
such that
\begin{equation*}
\ell_A\bigl(I+r_A(I^*)\bigr)=\ell_A(I)\cap\ell_A\bigl(r_A(I^*)\bigr)
                            =\ell_A(I)\cap r_A\bigl(\ell_A(I)\bigr)^*=0.
\end{equation*}
It therefore suffices to show that $J$ is closed and that $f$ extends to a bounded
right $A$-module homomorphism on~$J$; Proposition~\ref{prop:qmaxs-sigf-aw} will then
yield the result.

Suppose that $\seqn x$ and $\seqn y$ are sequences in $I$ and in $r_A(I^*)$, respectively
such that $(x_n+y_n)_{n\in\NN}$ converges to $z\in A$. Take an approximate identity
$(e_\alpha)$ in $I^*\cap I$ and let $\eps>0$. There is $n_0\in\NN$ such that
\begin{equation*}
\|x_n+y_n-x_m-y_m\|<\eps\qquad(n,m\geq n_0)
\end{equation*}
and for fixed $n,m\geq n_0$ choose $\alpha$ such that
\begin{equation*}
\|x_n-e_\alpha x_n\|<\eps\quad\text{and}\quad\|x_m-e_\alpha x_m\|<\eps.
\end{equation*}
Since $e_\alpha\,\ell_A(I)^*=0$, it follows that
\begin{equation*}
\begin{split}
\|x_n-x_m\| &\leq \|x_n-e_\alpha x_n\|+\|e_\alpha x_n-e_\alpha x_m\|
                  +\|e_\alpha x_m-x_m\|\\
            &= \|x_n-e_\alpha x_n\|+\|e_\alpha(x_n+y_n-x_m-y_m)\|
                  +\|e_\alpha x_m-x_m\|<3\,\eps
\end{split}
\end{equation*}
for all $n,m\geq n_0$. Hence the sequence $\seqn x$ converges to some $x\in I$
and therefore $\seqn y$ converges to $y=z-x$ in $r_A(I^*)$. Consequently, $J$
is closed.

Define $\tilde f\colon J\to A$ by $\tilde f(x+y)=f(x)$ whenever $x\in I$ and
$y\in\ell_A(I)^*$. The argument in the previous paragraph shows that $\tilde f$
is a bounded right $A$-module homomorphism extending~$f$. Since
\begin{equation*}
\begin{split}
\|\tilde f(x+y)\|^2 &= \|f(x)\|^2\leq\|f\|^2\,\|x\|^2\leq\|f\|^2\,\|x^*x+y^*y\|\\
                    &= \|f\|^2\,\|(x+y)^*(x+y)\|=\|f\|^2\,\|x+y\|^2,
\end{split}
\end{equation*}
it follows that $\|\tilde f\|=\|f\|$.
\end{proof}

Since there exist \sigfinite/ non-injective \AW*s, Proposition~\ref{prop:sigf-aw-selfinj}
shows that there is no Baer criterion \cite[\S\,3, 3.7]{LamMR} for operator modules.

\smallskip
An \AW* $A$ is said to be {\it normal\/} if each increasing net $(p_\alpha)$
of projections in $A$ has a supremum in~$A_{sa}$. Note that, in this case,
$\sup_\alpha p_\alpha=\bigvee_\alpha p_\alpha$, where the latter is the supremum
in~$P(A)$.

We shall make use of the following result, compare \cite{AGoldnfv}, \cite{Wrightezt};
a concise proof can be found in~\cite{Saito}.

\begin{theorem}\label{thm:sigfinite-is-normal}
Every \sigfinite/ \AW* is normal.
\end{theorem}

We proceed to apply this to our study of essential one-sided ideals started in
Section~\ref{sect:ess-right-ideals}.

\begin{proposition}\label{prop:awstar-normal}
Let\/ $A$ be an \AW* with the property that every closed right ideal\/ $I$ of\/ $A$
with\/ $\ell_A(I)=0$ is $\Oone$-essential. Then\/ $A$ must be normal.
\end{proposition}

\begin{proof}
Let $(p_\alpha)$ be an increasing net of projections in $A$ and let
$p=\bigvee_\alpha p_\alpha$. Then
\[
I=\ol{\bigcup_\alpha p_\alpha A}+(1-p)A
\]
is a closed right ideal of $A$ with $\ell_A(I)=0$. By assumption, $\ell_{I(A)}(I)=0$
and hence $\ell_{\bar A}(I)=0$ (where $\Abar$ denotes the regular monotone completion
of~$A$). It follows that, in $P(\Abar)$, $\bigvee p_\alpha +1-p=1$ whence
$\bigvee p_\alpha=p$. Since $\Abar$ is monotone complete, this implies that
$p=\bigvee p_\alpha={\Abar}_{sa}\text{-}\sup_\alpha p_\alpha$ in ${\Abar}_{sa}$
and hence $p=\sup_\alpha p_\alpha$ in $A_{sa}$.
\end{proof}

\begin{corollary}\label{cor:ess-right-ideals-in-sigfinite-aw}
A closed right ideal in a \sigfinite/ \AW*\/ $A$ is $\Oone$-essential if and only if it has
zero left annihilator in~$A$.
\end{corollary}

\begin{proof}
Let $I$ be a closed right ideal in $A$ such that $\ell_A(I)=0$.
We have seen that there is $a\in I$, $a\geq0$ such that
$\ell_A(\ol{aA})=\ell_A(\{a\})=0$. Take a countable approximate identity of $\ol{aAa}$
consisting of projections~$p_n$. Then $\bigvee_n p_n=1$ so, by
Theorem~\ref{thm:sigfinite-is-normal}, $\sup_n p_n=1$ in $A_{sa}$. Since the
embedding $A\hookrightarrow\Abar$ is $\sup$-preserving,
$\bigvee p_n={\Abar}_{sa}\text{-}\sup_\alpha p_\alpha$ in $\Abar$. It follows that
$p_{\,\ol{aAa}}^{}={\Abar}_{sa}\text{-}\sup_\alpha p_\alpha=1$ (see page~\pageref{p:B})
and so $\ell_{I(A)}(I)=0$, that is, $I$ is essential in~$\Oone$.
\end{proof}

%

\begin{proposition}\label{prop:small_awstar}
Let\/ $A$ be a separable \C*. Then there is a smallest {\sl AW*}-subalgebra of\/
$I(A)$ containing~$A$. This is a \sigfinite/ {\sl AW*}-subalgebra of~$\Abar$.
\end{proposition}

\begin{proof}
Let $\hat A$ be the regular monotone $\sigma$-completion of~$A$, which is
a \sigfinite/ \C* \cite{Wrightsvsa}. In the case considered here, $\hat A=\Abar$.
Let $B=\bigcap A'$, where $A'\subseteq\Abar$ is an \AW* containing~$A$. Any
such algebra $A'$ is \sigfinite/ and hence normal, by
Theorem~\ref{thm:sigfinite-is-normal}. It follows that $A'$ is an
{\sl AW*}-subalgebra of both $\Abar$ and $I(A)$, thus suprema of any subset of projections
in $A'$ are the same in $A'$, $\Abar$ and~$I(A)$. In particular,
the left support projections $s_{\ell,A'}(b)$ and $s_{\ell,I(A)}(b)$ agree for
all $b\in A'$. It follows that $B$ is an {\sl AW*}-subalgebra of~$\Abar$.

It remains to verify that $B\subseteq D$ for every {\sl C*}-subalgebra $D$ of
$I(A)$ containing $A$ which is an \AW*. As $A$ is separable, it acts faithfully
on a separable Hilbert space $H$ and thus $I(A)$ is completely isometric to
an operator system in $B(H)$. Hence $I(A)$ possesses a faithful state and
thus is \sigfinite/; it follows that $D$ is
\sigfinite/ and hence normal, by Theorem~\ref{thm:sigfinite-is-normal}.
Since $I(A)=I(D)$, we conclude that
\[
s_{\ell,D}(b)=s_{\ell,I(D)}(b)=s_{\ell,I(A)}(b)=s_{\ell,\bar A}(b)
\]
so that $\Abar\cap D$ is an \AW*. Therefore $B\subseteq\Abar\cap D\subseteq D$.
\end{proof}

We will now obtain a neat description of algebraically essential one-sided
ideals of finite \AW*s. This can be obtained from results by Utumi and Berberian
but we provide a direct argument making use of the fact that the projection lattice
$P(A)$ of a finite \AW* $A$ is continuous.

\begin{proposition}\label{prop:ess-rideals-in-finite-aw}
An \AW*\/ $A$ is finite if and only if every (not necessarily closed) right ideal\/
$I$ of\/ $A$ such that\/ $\ell_A(I)=0$ is algebraically essential in~$A$.
\end{proposition}

\begin{proof}
Suppose that $A$ is finite. Since $P(A)$ is a continuous lattice, for every orthogonal
family $(e_\alpha)$ in $P(A)$ the right ideal $\bigoplus_\alpha e_\alpha A$ is algebraically
essential in $(\bigvee_\alpha e_\alpha)A$. Let $I$ be a non-zero right ideal of~$A$.
Take a maximal orthogonal family of projections $e_\alpha\in I$ and put
$e=\bigvee_\alpha e_\alpha$. As
$\bigoplus_\alpha e_\alpha A+(1-e)A$ is algebraically essential in~$A$,
$\bigoplus_\alpha e_\alpha A+(1-e)A\cap I$ is algebraically essential in~$I$.
Since $(1-e)A\cap I=0$, it follows that $\bigoplus_\alpha e_\alpha A$ is
algebraically essential in~$I$. Therefore, $\ell_A(I)=A(1-e)$ and thus
$\bigoplus_\alpha e_\alpha A\subseteq I\subseteq eA=r_A(\ell_A(I))$.
As a result, $I$ is algebraically essential in $r_A(\ell_A(I))$, so in particular,
if $\ell_A(I)=0$ then $I$ is algebraically essential in~$A$.

Now suppose that $A$ is not finite; then $P(A)$ is not continuous.
Let $(e_\alpha)$ be an increasing net of projections and $q\in P(A)$ be such that
$\bigl(\bigvee_\alpha e_\alpha\bigr)\wedge q\ne\bigvee_\alpha(e_\alpha\wedge q)$.
Then $q_0=\bigl(\bigvee_\alpha e_\alpha\bigr)\wedge q-\bigvee_\alpha(e_\alpha\wedge q)$
is a non-zero projection such that, for all~$\alpha$, $e_\alpha\wedge q_0=0$.
Hence $\bigl(\bigcup_\alpha e_\alpha A\bigr)\cap q_0A=0$
which entails that the right ideal $\bigcup_\alpha e_\alpha A$ is not algebraically essential
in $\bigl(\bigvee_\alpha e_\alpha\bigr)A$.
\end{proof}

\begin{remark}\label{rem:ess_right_ideals_in_bh}
In contrast to the above proposition, if $H$ is a separable infinite dimensional
Hilbert space, there are increasing sequences $(e_n)$ of projections in $B(H)$
such that $\bigcup_n e_nB(H)$ is not algebraically essential in $B(H)$ but its closure
is. In fact, every closed right ideal in $B(H)$ with zero left annihilator is
(algebraically) essential in this case.
\end{remark}

In the finite case, we can also determine the maximal \C* of quotients.

\begin{theorem}\label{thm:qmax-finite-aw}
For every finite \AW*\/ $A$ we have\/ $\qmax=A$.
\end{theorem}

The proof is essentially contained in the literature, cf.~\cite{maxBer},
but we will guide the reader on how to put the pieces together. 
It depends on the construction and the properties of the regular ring of a finite \AW*. 
In the case where $A$ is a finite von Neumann algebra this ring is the algebra of unbounded
operators affiliated to $A$, a construction that goes back to the
seminal work of Murray and von Neumann~\cite{MvN}. In the general
case of a finite \AW*, the regular ring $R$ of $A$ was constructed
by Berberian and has the following two crucial properties (see~\cite{regBer} 
and~\cite{matrixBer}):
\begin{enumerate}
\item[(a)] $R$ is a *-regular ring containing~$A$;
\smallskip
\item[(b)] for every $n\ge1$, the \AW* $M_n(A)$ contains all the partial
isometries of the *-regular ring~$M_n(R)$.
\end{enumerate}
Indeed, the construction of the *-regular ring $R$ was a fundamental tool in the proof 
of the fact that matrix algebras over \AW*s are also \AW*s. 

\begin{proof}
Note that the projection lattice of
$M_2(R)$ is the same as the projection lattice of $M_2(A)$. Since
$M_2(A)$ is a finite \AW*, its projection lattice is a continuous
geometry and thus, by \cite[Corollary~13.19]{vnrr}, $R$ is
(algebraically) right and left self-injective. By construction,
for each element $a\in R$, there is an increasing sequence
$(e_n)_{n\in\NN}$ of projections in $A$ with supremum $1$ such that 
$ae_n\in A$ for every~$n$. Consequently, the right ideal 
$I=\bigcup_n e_nA$ satisfies $\ell_A(I)=0$ and hence is algebraically
essential, by Proposition~\ref{prop:ess-rideals-in-finite-aw}.
Using this it is fairly easy to show that $R\subseteq\qmaxs$. 
Since $R$ is self-injective, it follows that $R=\qmaxs$. Indeed, if
$g\colon I\to A$ is a right module homomorphism from a (not
necessarily closed) essential right ideal $I$ of $A$ into $A$,
then it can be extended to a right $R$-module homomorphism
$\widetilde{g}\colon IR\to R$ so that, by the injectivity of $R_R$, we
get an element $y\in R$ with the property that $\widetilde{g}(z)=yz$ for all
$z\in IR$; thus $\qmaxs\subseteq R$. 

Suppose that $a\in\qmaxsb$; then, by the above discussion and 
Remark~\ref{rem:ess-right-ideals-vs-algebras}, $a$ has a representative
$(f,g,I)$, with $I$ as above, and, moreover, $\|g(e_n)\|\le K$ for all~$n$. 
It follows from \cite[Theorem~5.1]{regBer} that $a\in A$ and thus $\qmax =A$, as desired.
\end{proof}

In a recent paper Argerami and Farenick determined those separable \C*s with the property
that $\mloc$ or $I(A)$ is a von Neumann algebra~\cite{ArgFar1}. It turns out that, for the
same class of \C*s, $\qmax$ is a von Neumann algebra.

\begin{proposition}\label{prop:qmax_vN}
Let\/ $A$ be a separable \C*. Then\/ $\qmax$ is a von Neumann algebra if and only
if\/ $A$ contains a minimal essential ideal consisting of compact elements.
\end{proposition}

Recall that an element $a\in A$ is {\it compact\/} if the two-sided multiplication
operator $x\mapsto axa$ is compact on~$A$. The set of all compact elements in $A$
is a closed ideal which is either zero or a direct sum of elementary \C*s, that is,
copies of the algebra of compact operators. See, e.g., \cite[Proposition~1.2.30]{AMMzth}.

\begin{proof}
Suppose that $K=\bigoplus_n K(H_n)$ is a minimal essential ideal of compact elements
in~$A$. By \cite[Lemma~1.2.21 and Proposition~2.3.6]{AMMzth},
\[
\mloc=M(K)=M\bigl(\bigoplus_n K(H_n)\bigr)=\prod_n M(K(H_n))=\prod_n B(H_n)
\]
and thus $\mloc$ is an atomic type~I von Neumann algebra, in particular injective.
By Proposition~3.5 and Theorem~3.8 it follows that $\mloc=\qmax=I(A)$ is a von Neumann algebra.

\medbreak\goodbreak

Conversely, suppose that $\qmax$ is a von Neumann algebra.
Propositions~\ref{prop:qmax_into_abar} and~\ref{prop:invariant-injective} together with
the hypothesis imply that $\Abar=\ol{\qmax}=\qmax$ is a von Neumann algebra.
Since the pure state space of $\Abar$ is separable~\cite{Wrightsvsc}, it follows
from \cite[Corollary~7]{Wrightsvsb} that $\Abar=\prod_n B(H_n)$ for a sequence
$(H_n)_{n\in\NN}$ of separable Hilbert spaces. Section~3 in \cite{Hameta} finally
entails that $K=\bigoplus_n K(H_n)$ is a minimal essential ideal in~$A$.
\end{proof}

\section{An Example}\label{sect:example}

In this section we shall study the details of an example in order to get a feel
for the properties and the behaviour of the local multiplier algebra and the maximal
\C* of quotients. This will lead us to a new \C* $A$ with the property that
$\Mloc\mloc\ne\mloc$.

We begin by recalling that, for a commutative unital \C* $A=C(X)$, we have
$\mloc=\qmax=B(X)$, the \C* of bounded complex-valued Borel
functions on $X$ modulo the ideal of those functions vanishing outside a set
of first category \cite[Proposition~3.4.5]{AMMzth}. Since $B(X)$ is an injective
\AW*, it follows that $\mloc=\qmax=\Abar=I(A)$. It is well known that, if $X$
contains a dense subset of first category, then $B(X)$ is not a von Neumann
algebra \cite[Proposition~III.1.26 and Theorem~III.1.17]{TAK}; if $X$ is
second countable, then $B(X)$ is \sigfinite/.

Our example is built from the commutative case and a non-commutative \C* in the
most basic manner. Let $X$ be a compact Hausdorff space and let $H$ be an infinite-dimensional
Hilbert space. Let $A=C(X)\otimes B(H)=C(X,B(H))$ be the {\sl C*}-tensor product
of $C(X)$ and~$B(H)$.

We start by computing the local multiplier algebra of~$A$.

\begin{proposition}\label{prop:2-sidedideals}
Every closed essential ideal of\/ $A$ contains an ideal of the form\/
$C_0(U,K(H))$, where\/ $\,U$ is an open dense subset of\/ $X$.
\end{proposition}

\begin{proof}
Let $I$ be a closed essential ideal of~$A$. Then
\[
Y=\{t\in X\mid f(t)=0\,\text{ for all }f\in I\}
\]
is closed in $X$ and $U:=X\setminus Y$ is open. Since $I$ is essential, $U$ must be
dense in~$X$. We want to show that $C_0(U,K(H))\subseteq I$. Since the algebraic
tensor product $C_0(U)\odot K(H)$ is dense in $C_0(U,K(H))$,
it is enough to show that $C_0(U)\odot K(H)$ is contained in~$I$. As
the linear combinations of elements $f\otimes e$, where $f\in C_0(U)$ and $e$ is a
one-dimensional projection, are dense in $C_0(U)\odot K(H)$, it suffices to prove
that $f\otimes e\in I$ for all $f\in C_0(U)$ and any one-dimensional projection~$e$.
We identify $f\otimes e$ with the function $t\mapsto f(t)e$ in $C_0(U,K(H))$,
written as~$fe$.

The closed ideal
\[
J=\{f\in C(X)\mid fe\in I\}
\]
is of the form $J=C_0(V)$ for some open subset $V\subseteq X$.
Clearly $V\subseteq U$ so we need to show that $V=U$. Take $t_0\in U$.
There exists $g\in I$ such that $g(t_0)\ne 0$. Consequently there exist $a,b\in F(H)$,
the finite-rank operators in $B(H)$, such that
$ag(t_0)b=e$. Identifying $a$ and $b$ with the constant functions
$t\mapsto a$ and $t\mapsto b$, respectively, we get $agb\in I$
because $g\in I$ and $(agb)(t_0)=ag(t_0)b=e$. Moreover $(ea)g(be)\in I$ and,
for all $t\in X$, we have
\[
eag(t)be=h(t)e
\]
for some $h(t)\in\CC$. Since $t\mapsto h(t)e$ is continuous, the function $h$ must 
be continuous. Therefore $h\in J$ and $h(t_0)\ne 0$, showing that $V=U$. 
We conclude that $C_0(U,K(H))\subseteq I$.
\end{proof}

\begin{remark}\label{rem:ideals-in-min-tensorp}
A closed ideal in the minimal tensor product $A_1\otimes_{\text{\rm min}}A_2$
of two unital \C*s $A_1$ and $A_2$ contains a closed ideal of the form
$I_1\otimes I_2$, where $I_j\subseteq A_j$, $j=1,2$ are closed ideals;
see \cite[Exercise~IV.4.3]{TAK}. The above is an easy direct argument
in the case considered here.
\end{remark}

For a locally compact Hausdorff space $U$, the multiplier algebra of $C_0(U,K(H))$
is given by $M\bigl(C_0(U,K(H))\bigr)=C_b(U,B(H)_\beta)$,
where $\beta$ stand for the strict topology \cite[Corollary~3.5]{APT}.
On $B(H)$, this agrees with the strong*-topology.
(We note that, on bounded subsets of $B(H)$, it also coincides with
the Arens--Mackey topology $\tau(B(H),B(H)_*)$, the finest
locally convex topology making all functionals in $B(H)_*$ continuous
\cite[Theorem~III.5.7]{TAK}.) Combining this fact with
Proposition~\ref{prop:2-sidedideals}, we obtain the description of $\mloc$.

\begin{corollary}\label{cor:locmult-cxbh}
Let\/ $A=C(X,B(H))$. Then
\[
\mloc=\Dirlim{U\in\mathfrak D}\;C_b(U,B(H)_\beta),
\]
where\/ $\mathfrak D$ is the filter of dense open subsets of\/ $X$.
\end{corollary}

In order to identify the maximal \C* of quotients we need more information
on one-sided ideals.

\begin{lemma}\label{lem:openxi}
Let\/ $I$ be an essential closed right ideal of~$A$. For each unit vector $\xi\in H$
there exists an open dense subset\/ $U_{\xi}$ of\/ $X$ such that
$C_0(U_{\xi})e_{\xi}A\subseteq I$, where\/ $e_\xi$ denotes the orthogonal projection
onto the subspace spanned by\/~$\xi$.
\end{lemma}

\begin{proof}
The closed ideal
\[
J=\{f\in C(X)\mid fe_{\xi}\in I\}
\]
is of the form $J=C_0(U_{\xi})$ for an open subset $U_{\xi}\subseteq X$. We verify that
$U_{\xi}$ is dense in~$X$. If not, then $X\setminus\overline{U_{\xi}}$ is open and non-empty
so there is $g\in C(X)$ such that $g\ne 0$ and $g(\overline{U_{\xi}})=0$. Since $I$ is
(algebraically) essential, there is $h\in A$ such that $ge_{\xi}A\cap I\ni ge_\xi h\ne0$.
Observe that
\[
0\ne ge_{\xi}hh^*e_{\xi}=gke_{\xi}
\]
for some $k\in C(X)$. As this implies that, at the same time, $0\ne gk\in J=C_0(U_{\xi})$ and
$(gk)(\ol{U_{\xi}})=0$, we obtain a contradiction. This shows the result.
\end{proof}

We now show that $\qmax$ agrees with $\mloc$ under an additional assumption on the
topological space~$X$.

\begin{proposition}\label{prop:mloc-is-qmax-if-twostar}
Let\/ $A=C(X,B(H))$ for a compact Hausdorff space\/ $X$ with the property that every
subset of first category is rare and a separable Hilbert space\/~$H$. Then
every essential closed right ideal of\/ $A$ contains an essential closed ideal.
Consequently\/ $\mloc=\qmax$.
\end{proposition}

\begin{proof}
Take an orthonormal basis $\{\xi_n\mid n\in\NN\}$ of~$H$ and let $I\in\Icer$.
By Lemma~\ref{lem:openxi}, we can find open dense subsets $U_n$ of $X$ such that
$C_0(U_n)e_{\xi_n}A\subseteq I$ for all~$n$. Let $V$ be
the interior of $\bigcap_{n=1}^{\infty}U_n$. By the hypothesis on~$X$,
$V$ is an open dense subset of $X$. We will show that
$C_0(V,K(H))\subseteq I$. As in the proof of Proposition~\ref{prop:2-sidedideals},
it suffices to show that $f\otimes e\in I$ for each $f\in C_0(V)$ and each one-dimensional
projection $e\in B(H)$. Fix a one-dimensional projection $e$ and consider the closed ideal
\[
J=\{f\in C(X)\mid fe\in I\},
\]
which is of the form $J=C_0(W)$ for some open subset $W\subseteq X$. We want
to show that $V\subseteq W$. Observe that
$\ol{\sum_{n=1}^{\infty} C_0(V)e_{\xi_n}A}\subseteq I$. Take $t_0\in V$ and $f\in C_0(V)$
such that $f(t_0)=1$. Since $\ol{\sum_{n=1}^{\infty} e_{\xi_n}B(H)}=K(H)$, we can find
$a_i\in\sum_{n=1}^{\infty}e_{\xi_n}B(H)$, $i\in\NN$ such that $e=\lim_ia_i$.
As $fe=\lim_i fa_i$ in $A$, it follows that
\[
fe\in\ol{\textstyle{\sum_{n=1}^{\infty}C_0(V)e_{\xi_n}A}}\subseteq I.
\]
Thus $f\in J$ and $f(t_0)=1$ implying that $t_0\in W$, and we conclude that
$V\subseteq W$ as desired.
\end{proof}

The topological condition in the above proposition arises naturally in the context
of commutative \AW*s. Let $X$ be a Stonean space (that is, $X$ is a compact Hausdorff
space such that the closure of every open subset is open). Then $X$ is the disjoint
union of two open and closed subsets $X_i$, $i=1,2$, where $X_1$ contains a dense
subset of first category and every subset of $X_2$ of first category is rare (that is,
nowhere dense). The set $X_2$ can further be decomposed into a disjoint union of
closed open subsets $X_{21}$ and $X_{22}$, where $X_{21}$ is hyper-Stonean (that is,
the spectrum of a commutative von Neumann algebra) and the support of every measure
on $X_{22}$ is rare. See \cite[Theorem~1.17]{TAK}. An example of a subset $X_1$
is the spectrum of $\Mloc{C[0,1]}$, see the remarks at the beginning of this section.

Letting $\mathfrak T$ be the set of dense $G_{\delta}$-subsets of the
compact space $X$ we have for the local multiplier algebra of $A=C(X,B(H))$:
\begin{equation}
\mloc=\Dirlim{U\in\mathfrak D}\;C_b(U,B(H)_\beta)\subseteq\Dirlim{T\in\mathfrak T}\;
C_b(T,B(H)_\beta).
\end{equation}
Observe that the latter direct limit of \C*s is indeed an algebraic direct limit, since
$\mathfrak T$ is closed under countable intersections.

In our next theorem we shall characterise the right-most algebra in~(1) above
as the injective envelope of $A$, provided $X$ is Stonean. To this end we recall
some results by Hamana~\cite{Hametb}, \cite{Hametc}.

Suppose that $A=C(X)\otimes B(H)$, where $X$ is a Stonean space. Then Hamana
introduces the monotone complete tensor product $C(X)\botimes B(H)$ and shows
that $\Abar=I(A)=C(X)\botimes B(H)$. Moreover, $I(A)$ is a von Neumann algebra
precisely when $X$ is hyper-Stonean. See \cite[Theorem~4.2 and Corollary~4.11]{Hametb}.
Denoting by $B(H)_\sigma$ the space $B(H)$ endowed with the $\sigma$-weak topology
$\sigma(B(H),B(H)_*)$ we have $I(A)=C(X,B(H)_\sigma)\,$ by \cite[Theorem~1.3]{Hametc}.

\begin{theorem}\label{thm:ia-cxbh}
Let\/ $A=C(X)\otimes B(H)$ for a Stonean space\/ $X$ and a separable Hilbert space\/~$H$.
Then\/ $I(A)=\Dirlim{T\in\mathfrak T}\;C_b(T,B(H)_\beta)$.
\end{theorem}

\begin{proof}
We first show the inclusion
$\Dirlim{T\in\mathfrak T}\;C_b(T,B(H)_\beta)\subseteq C(X,B(H)_\sigma)$.
Let $f\in C_b(T,B(H)_\beta)\subseteq C_b(T,B(H)_\sigma)$ for some $T\in\mathfrak T$.
Since $X$ is Stonean and bounded subsets of $B(H)$ are relatively $\sigma$-weakly compact,
$f$ can be uniquely extended to a function in $C(X,B(H)_\sigma)$.

For the reverse inclusion, let $\{\xi_n\mid n\in\NN\}$ be an orthonormal basis
of~$H$ and take $f\in I(A)_+$. Putting $\xi(t)=f(t)\eta$ for a fixed vector $\eta\in H$
and all $t\in X$, the continuity assumption on $f$ yields that
$t\mapsto(\xi(t)\mid\xi_n)$ is continuous on $X$ for each~$n$.
Since $\|\xi(t)\|^2=\sum_{n=1}^\infty|(\xi(t)\mid\xi_n)|^2$, it follows that
$t\mapsto\|\xi(t)\|^2$ is lower semicontinuous; thus it is continuous on a dense
$G_{\delta}$-subset $T$ of~$X$. From the identity
\[
\|\xi(s)-\xi(t)\|^2=\|\xi(s)\|^2+\|\xi(t)\|^2-2\,\textit{Re}\,(\xi(s)\mid\xi(t))
\qquad(s,t\in X)
\]
we infer that $\xi\in C_b(T,H)$. Applying this argument to each element in
a countable dense subset $S$ of $H$,
we obtain a countable family $\{T_k\mid k\in\NN\}$ of dense $G_\delta$-subsets
of $X$ such that $t\mapsto f(t)\eta_k$ is continuous on~$T_k$. Letting
$T=\bigcap_k T_k$ we obtain a dense $G_\delta$-subset on which all these functions
are simultaneously continuous.

Let $\eta\in H$. We claim that $t\mapsto f(t)\eta$ is continuous on~$T$.
Take $t_0\in T$ and let $\eps>0$. Choose
$\eta'\in S$ such that $\|\eta-\eta'\|<\epsilon$. Since
$t\mapsto f(t)\eta'$ is continuous on $T$, there is an open neighbourhood $V$ of $t_0$
in $T$ such that $\|f(t)\eta'-f(t_0)\eta'\|<\epsilon$ for all $t\in V$. For
$t\in V$ we have
\begin{equation*}
\begin{split}
\|f(t)\eta-f(t_0)\eta\|
&\le \|f(t)\eta-f(t)\eta'\|+\|f(t)\eta'-f(t_0)\eta'\|+\|f(t_0)\eta'-f(t_0)\eta\|\\
&< (2\,\|f\|+1)\,\epsilon.
\end{split}
\end{equation*}
Consequently, $f\in C_b(T,B(H)_\beta)$ as claimed.
\end{proof}

\begin{corollary}\label{cor:mloc-is-ia-if-twostar}
Let\/ $A=C(X,B(H))$ for a Stonean space\/ $X$ with the property that every
subset of first category is rare and a separable Hilbert space\/~$H$. Then\/
$\mloc=I(A)$. Moreover, $\mloc$ is a von Neumann algebra if and only if\/ $X$
is hyper-Stonean.
\end{corollary}

\begin{proof}
Since the hypothesis on $X$ implies that every dense $G_\delta$-subset of $X$ contains a
dense open subset, the statement follows immediately from Theorem~\ref{thm:ia-cxbh}
and Corollary~\ref{cor:locmult-cxbh}.
\end{proof}

\begin{remark}
The proof of Theorem~\ref{thm:ia-cxbh} is inspired by the argument in \cite{Hametc},
p.~291 where it is shown that every subset of first category in $X$ is rare if and only
if every projection in $\Abar$ is open. This, in turn is equivalent to the
requirement that each positive element in $\Abar_+$ is the supremum of an increasing
net in $A_+$, by \cite[Lemma~1.7]{Hameta}.
\end{remark}

In a \C* of the form $A=C(X,B(H))$, where $X$ is a Stonean space and $H$ is a
separable Hilbert space, we find $\Oone$-essential right ideals which
are not algebraically essential.

\begin{example}\label{exam:o1ess-not-alg-ess}
Suppose that $X$ contains a decreasing sequence $(U_n)_{n\in\NN}$
of open dense subsets such that $\bigcap_{n=1}^{\infty} U_n$ has empty interior.
Fix an orthonormal basis $\{\xi_n\mid n\in\NN\}$ of $H$ and set
$I=\ol{\sum_{n=1}^\infty C_0(U_n)e_{\xi_n}A}$,
which is a closed right ideal of~$A$. Then $I$ is not
algebraically essential in~$A$.
(This part does not use that $X$ is Stonean, only that $X$ is a compact
Hausdorff space.) Note that every element of $I$ is of the form
\[
\sum_{n=1}^{\infty}\xi_n\otimes\eta_n,
\]
where $\eta_n\in C_0(U_n,H)$ for all $n$ and the convergence is
uniform. Put $\xi=\sum_{n=1}^\infty 2^{-n}\xi_n$ and
$f=\xi\otimes\xi\in A$. If $I$ were algebraically essential then
$fA\cap I\ne 0$; thus we can write
\[
0\ne\xi\otimes\eta=\sum_{n=1}^{\infty}\xi_n\otimes\eta_n,
\]
where $\eta\in C(X,H)$
and $\eta _n\in C_0(U_n,H)$ for all~$n$. Since $\eta$ is
continuous, there is a non-empty open subset $G$ of $X$ such that
$\eta(t)\ne0$ for all $t\in G$. By hypothesis $G\nsubseteq\bigcap_{n=1}^{\infty}U_n$,
so that we can choose $t_0\in G$ and $n_0\in\NN$ such that
$t_0\notin U_n$ for all $n\ge n_0$. We obtain
\[
0\ne\xi\otimes\eta(t_0)=\sum_{n=1}^\infty\xi_n\otimes\eta_n(t_0)
                         =\sum_{n=1}^{n_0}\xi_n\otimes\eta_n(t_0),
\]
and so $\|\eta(t_0)\|^2\,\xi=\sum_{n=1}^{n_0}(\eta(t_0)\mid\eta_n(t_0))\,\xi_n$,
a contradiction to the definition of~$\xi$.

\smallskip
To see that $I$ is $\Oone$-essential, it suffices to show that
$\ell_{I(A)}(I)=0$; see Theorem~\ref{thm:ess-right-ideals}. Let
$f$ be a non-zero element of~$I(A)$. By Theorem~\ref{thm:ia-cxbh},
there is a dense $G_{\delta}$-subset $T$ of $X$ such that
$f\in C_b(T,B(H)_{\tau})$. On multiplying $f$ by a function of
the form $\xi\otimes\xi$, for suitable $\xi\in H$, on the left, we get
$0\ne\xi\otimes f(t)^*\xi $ with $t\mapsto\xi(t):=f(t)^*\xi$
norm-continuous on~$T$. There exist a natural number $m$ and a
non-empty subset $V$ of $X$ such that $(\xi(t)\mid\xi_m)\ne 0$ for
all $t\in V\cap T$. Since $T$ and $U_m$ are dense, we obtain
$U_m\cap V\cap T\ne \emptyset $. Take $t_0\in U_m\cap V\cap T$ and
$h\in C_0(U_m)_+$ such that $h(t_0)\ne0$. Then
$h(\xi_m\otimes \xi_m)$ is in $I$ and
$(\xi \otimes \xi(t))\,h(\xi_m\otimes \xi_m)$ is non-zero because
\[
(\xi\otimes\xi(t_0))\,h(t_0)(\xi_m\otimes \xi _m)
=\xi \otimes h(t_0)(\xi(t_0)\mid\xi_m)\,\xi_m\ne0.
\]
\end{example}

We now want to add some remarks on the Hilbert {\sl C*}-module structure of
the spaces discussed above. For the remainder of this section,
$B=C(X)$ for a Stonean space~$X$ and $A=B\otimes B(H)$, where $H$ is a
separable Hilbert space with orthonormal basis $\{\xi_n\mid n\in\NN\}$.

The space $E=C(X,H_w)$ is a faithful \AWm* over the \AW* $B$ in a canonical way
(in the sense of Kaplansky~\cite{Kapfvth}); the $B$-valued inner product is given by
\[
(\xi\mid\eta)_E^{}(t)=(\xi(t)\mid\eta(t))_H^{}
\]
for all~$t$ in a dense $G_\delta$-subset of~$X$. We shall normally drop the
subscripts `$E$' and `$H$', if no confusion can arise.
Note that the norm in $E$ coincides with the sup-norm and that $E$ is unitarily
equivalent with $\ell^2(B)$ via $\xi\mapsto(\xi\mid\xi_n)$ so that
$\{\xi_n\mid n\in\NN\}$ is an orthonormal basis of the \AWm* $E$ (we continue to identify
$\eta\in H$ with the constant function $t\mapsto\eta$, $t\in X$).
See \cite[Proposition~1.7]{Hametc}. Recall also that $C_b(U,H)$ embeds isometrically
into $E$ for every $U\in\mathfrak D$.

\begin{lemma}\label{lem:f-is-e}
Let\/ $G$ be the closure of the space\/
\[
\{\xi\in E\mid\rest{\xi}{U}\in C_b(U,H)\text{ for some dense open subset }\,U\subseteq X\}.
\]
Then\/ $G=E$.
\end{lemma}
\begin{proof}
We will establish the claim by showing that $G$ is an \AWsm* of $E$ with zero
orthogonal complement; the fact that every \AWsm* is orthogonally complemented
\cite[Theorem~3]{Kapfvth} will then complete the argument.

The continuity of the module operations immediately yields that $G$ is a closed
$B$-submodule. In order to show that $G$ is an \AWsm*, let
$\{\xi_i\}$ be a bounded subset of $G$ and let $\{e_i\}$ be an orthogonal family
of projections in~$B$ with supremum~$1$. We need to verify that $\sum_i e_i\xi_i\in G$.

To this end we can assume that, for each~$i$, $\rest{\xi_i}{U_i}\in C_b(U_i,H)$
for some open dense subset $U_i\subseteq X$. Suppose we can show that
$\sum_i e_i\xi_i\in G$ under this assumption. For the general case,
let $\eps>0$ and,
for each~$i$, choose $\xi_i'\in E$ such that $\rest{\xi_i'}{U_i}\in C_b(U_i,H)$
for some open dense subset $U_i\subseteq X$ and that $\|\xi_i-\xi_i'\|<\eps$.
Since
\[
\Bigl\|\sum_ie_i\xi_i-\sum_ie_i\xi_i'\Bigr\|=\sup_i\|e_i(\xi_i-\xi_i')\|\leq\eps
\]
and $\sum_i e_i\xi_i'\in G$, it follows that $\sum_i e_i\xi_i\in G$ as needed.

Making this additional assumption, let $V_i\subseteq X$ be mutually disjoint clopen
subsets of $X$ such that $e_i=\chi_{V_i}$ for each~$i$. Since $\sup_ie_i=1$,
$\bigcup_iV_i$ is dense in~$X$.
Since $\sum_i e_i\xi_i$ is the unique element in $E$ such that
$e_j\bigl(\sum_i e_i\xi_i\bigr)=e_j\xi_j$ for all~$j$; the element
$\rest{e_j\xi_j}{V_j\cap U_j}\in C_b(V_j\cap U_j,H)$ for all~$j$;
and the $V_j\cap U_j$'s are mutually disjoint open subsets,
it follows that $\rest{\sum_i e_i\xi_i}{U}\in C_b(U,H)$,
where $U=\bigcup_j(V_j\cap U_j)$ is open and dense.
As a result, $\sum_i e_i\xi_i\in G$.

We finally show that the orthogonal complement of $G$ is zero.
Suppose that $\xi=\sum_{n=1}^\infty(\xi\mid\xi_n)\,\xi_n\in E$ is non-zero;
say, $(\xi\mid\xi_m)\,\xi_m\ne0$. As
$\xi'=(\xi\mid\xi_m)\,\xi_m\in C(X,H)\subseteq G$ and
$(\xi\mid\xi')=|(\xi\mid\xi_m)|^2\ne0$, we conclude that no non-zero element
in $E$ is orthogonal to~$G$.
\end{proof}

Let $\LB(E)$ denote the \AW* of bounded (adjointable) module homomorphisms of~$E$.
By \cite{Hametc}, $\LB(E)=I(A)$. Let $\KB(E)$ denote the so-called compact elements
in $\LB(E)$, that is, the closed linear span of elements of the form
$\theta_{\xi,\eta}$, $\xi,\eta\in E$~\cite{Lance}. Since
$\|\theta_{\xi,\eta}\|\leq\|\xi\|\,\|\eta\|$ and
\[
\|\theta_{\xi,\eta}-\theta_{\xi',\eta'}\|\leq\|\xi-\xi'\|\,\|\eta\|+\|\xi'\|\,\|\eta-\eta'\|
\]
for all $\xi,\xi',\eta,\eta'\in E$, it follows from Lemma~\ref{lem:f-is-e}
that every $\theta_{\xi,\eta}$ is uniformly approximated by $\theta_{\xi',\eta'}$
with $\xi',\eta'\in C_b(U,H)$ for a dense open subset $U\subseteq X$.
Combining this with Corollary~\ref{cor:locmult-cxbh} we obtain the following result.

\begin{proposition}\label{prop:kb-in-mloc}
With the above notation and caveats we have\/ $\KB(E)\subseteq\mloc$.
\end{proposition}
\begin{proof}
The isomorphism between $\LB(E)$ and $I(A)=C(X,B(H)_\sigma)$ carries the operator
$\theta_{\xi,\eta}$, $\xi,\eta\in E$ onto the function $\xi\otimes\eta$ defined by
$(\xi\otimes\eta)(t)\,\zeta=(\zeta\mid\eta(t))\,\xi(t)$, $t\in X$, $\zeta\in H$.
If $\xi,\eta\in C_b(U,H)$ for a dense open subset $U\subseteq X$ then
$\xi\otimes\eta\in C_b(U,B(H)_\beta)$, since, for all $\zeta\in H$ and all $s,t\in U$,
\[
\|(\xi\otimes\eta)(t)\,\zeta-(\xi\otimes\eta)(s)\,\zeta\|
\leq
\|\zeta\|\,\|\eta(t)-\eta(s)\|\,\|\xi(t)\|+\|\zeta\|\,\|\eta(s)\|\,\|\xi(t)-\xi(s)\|
\]
and $(\xi\otimes\eta)^*=\eta\otimes\xi$.
Thus, identifying $\KB(E)$ with the closed linear span of
$\{\xi\otimes\eta\mid\xi,\eta\in E\}$ and using Lemma~\ref{lem:f-is-e}
together with Corollary~\ref{cor:locmult-cxbh}, it follows
that $\KB(E)\subseteq\mloc$.
\end{proof}

Lemma~\ref{lem:f-is-e} can be reformulated by stating that, as
Hilbert {\sl C*}-modules over $B$,
we have $C(X,H_w)=\Dirlim{U\in\mathfrak D}\,C_b(U,H)$.
Therefore, as \C*s, 
$\KB\bigl(C(X,H_w)\bigr)=\Dirlim{U\in\mathfrak D}\,\KB\bigl(C_b(U,H)\bigr)$.
Moreover,
\[
C_b(U,B(H)_\beta)=M\bigl(C_0(U,K(H))\bigr)\supseteq\KB\bigl(C_b(U,H)\bigr)\supseteq C_0(U,K(H)).
\]

We can depict the situation in the following commutative diagram.

\begin{equation*}
\specialsixdowndiag%
{\KB\bigl(C_b(U,H)\bigr)}{}{C_b(U,B(H)_\beta)}{\smalldirlim}%
{ }{\text{id}}{\KB\bigl(C(X,H_w)\bigr)}{\MLOC{C(X,B(H))}}{C(X,B(H)_\sigma)}
\end{equation*}

\medskip
Note that $\KB(E)$ is an essential ideal in $\mloc$; indeed, $I(A)=\LB(E)=M(\KB(E))$.
Therefore, $\Mloc\mloc=I(A)$.\label{p:mlocmloc}

\smallskip
In view of the above discussion
and Corollary~\ref{cor:mloc-is-ia-if-twostar} a natural question arises:
under what conditions is $\mloc\ne I(A)$? To answer this question we first need an
auxiliary result which undoubtedly is known; alas we do not have a reference.

\begin{lemma}\label{lem:element-in-mke}
With the above notation and caveats let\/ $\{\xi_n\mid n\in\NN\}$ be an orthonormal
basis of\/ $E$ and let $\{\eta_n\mid n\in\NN\}$ be an orthonormal family in~$E$.
Then\/ $a=\sum_{n=1}^\infty \theta_{\eta_n,\xi_n}$ defines an element
in~$\mathcal L_B(E)$.
\end{lemma}

\begin{proof}
To prove the claim of the lemma we need to show that both $a$ and
$a^*=\sum_{n=1}^\infty \theta_{\xi_n,\eta_n}$ define $B$-module maps from $E$ into~$E$.
Let $\zeta\in E$. Writing $\zeta=\sum_{m=1}^\infty c_m\xi_m$ with
$\sum_{m=1}^\infty c_m^*c_m^{}$ convergent in $B$ we obtain
$a\zeta=\sum_{n=1}^\infty c_n\eta_n\in E$, compare \cite[Lemma~9]{Kapfvth}.

Applying \cite[Theorem~1]{Kapfvth} to the orthogonal complement of the
{\sl AW*}-submodule of $E$ generated by $\{\eta_n\mid n\in\NN\}$ we obtain
an orthogonal family $\{e_i\}$ of projections in $B$ with supremum~$1$ such
that, for each~$i$, the family $\{e_i\eta_n\mid n\in\NN\}\cup\{\rho_{i,k}\mid k\in\NN\}$
is an orthonormal basis of the $e_iB$-module~$e_iE$, where the $\rho_{i,k}\in E$
are suitably chosen. For $\zeta\in E$ we now write
\[
\zeta=\sum_i e_i\zeta=\sum_i\sum_m b_{i,m}e_i\eta_m + \sum_i\sum_k c_{i,k}\rho_{i,k}
\]
with both $\sum_{m=1}^\infty e_i^{}b_{i,m}^*b_{i,m}^{}$ and
$\sum_{k=1}^\infty c_{i,k}^*c_{i,k}^{}$ bounded by $\|\zeta\|^2$ for all~$i$.
This yields that
\begin{equation*}
e_i \sum_{n=1}^\infty \theta_{\xi_n,\eta_n}(\zeta)
=\sum_{n=1}^\infty(e_i\zeta\mid e_i\eta_n)\xi_n
=\sum_{n=1}^\infty e_ib_{i,n}\xi_n\in e_iE
\end{equation*}
and thus
$a^*\zeta=\sum_{n=1}^\infty \theta_{\xi_n,\eta_n}(\zeta)
=\sum_i e_i\sum_{n=1}^\infty \theta_{\xi_n,\eta_n}(\zeta)\in E$.
\end{proof}

Note that the condition on the space $X$ in the next theorem implies that
$X$ contains a dense subset of first category.

\begin{theorem}\label{thm:mloc-not-ia}
Let\/ $X$ be a Stonean space with the property that there exists a
decreasing sequence\/ $(U_n)_{n\in\NN}$ of open dense subsets of\/
$X$ such that\/ $U_1=X$ and\/ $\bigcap_{n=1}^\infty U_n$ has empty
interior and, furthermore, that\/ every family of non-empty mutually disjoint
clopen subsets of\/ $X$ is countable. Let\/ $H$ be a separable (infinite-dimensional) 
Hilbert space, and let\/ $A=C(X,B(H))$.
Then\/ $\mloc$ is a proper \Cs* of\/~$I(A)$.
\end{theorem}

\begin{proof}
Without loss of generality, we can assume (using our hypothesis
that $C(X)$ is $\sigma$-finite) that each $U_n$ is a disjoint
union of countably many cozero sets of the form $\{t\in X\mid f_{\alpha}(t)>0\}$
for some $0\ne f_{\alpha}\in C(X)$ with $0\le f_\alpha\le1$.
It follows that $\chi_{U_n}$ equals the pointwise sum
$\sum_\alpha\bigl(f_\alpha+\sum_{i=1}^\infty(f_\alpha^{1/(i+1)}-f_\alpha^{1/i})\bigr)$ of
countable many continuous positive functions; see \cite[p.~291]{Hametc}.

For each $m\in\NN$, put $f_m=\sum_{i=1}^\infty \frac{1}{2^i}\chi_{U_{m-1+i}}$.
Then $f_m$ is bounded by~$1$ and lower-semicontinuous on $X$ and discontinuous at each point
of $X\setminus\bigcap_{n=1}^\infty U_n$.
By the above observation, there are non-zero positive functions $g_k^{(m)}\in C(X)$
such that $f_m(t)^2=\sum_{k=1}^\infty g_k^{(m)}(t)$ for all $t\in X$.
Let $\{\xi_n\mid n\in\NN\}$ be an orthonormal basis of~$H$ and put
$\zeta_m=\sum_{k=1}^\infty(g_k^{(m)})^{1/2}\xi_k\in E=C(X, H_w)$.
\smallskip
Then
$\|\zeta_m(t)\|^2=\sum_{k=1}^\infty g_k^{(m)}(t)=f_m(t)^2$ for all~$t$.
Moreover, as $f_m(t)=1$ for all $t$ in the dense $G_\delta$-subset $\bigcap_{n=1}^\infty U_n$,
we have $|\zeta_m|=1$ for all~$m$.

We claim that, for every $0<\eps<\frac{1}{12}$ and every open dense subset $U$
of $X$, there is $m\in\NN$
such that, for each $g\in C_b(U,H)$, the estimate $\|\zeta_m-g\|>\eps$ holds.

To prove the claim we shall show that, given $U$,
\begin{equation}\label{equ:bigg_than_eps}
\exists\;m\in\NN\ \forall\;h\in C_b(U)\ \exists\;t\in U\colon|f_m(t)-h(t)|>\eps.
\end{equation}

This indeed suffices, since, for $g\in C_b(U,H)$, the function $h\colon t\mapsto\|g(t)\|$
belongs to $C_b(U)$ and
\[
\|\zeta_m-g\|\geq\bigl|\|\zeta_m(t)\|-\|g(t)\|\bigr|=|f_m(t)-h(t)|\qquad(t\in X).
\]

Let $U\subseteq X$ be dense and open and let $\eps>0$.
Let $R_n=X\setminus U_n$, $n\in\NN$. By assumption, $\bigcup_{n=1}^\infty R_n$
is dense in $X$; thus there is a smallest $m\in\NN$ such that
$U\cap R_{m+1}\ne\emptyset$. In particular, $U\subseteq U_m$.

For $t_0\in U\cap R_{m+1}$ we have $f_m(t_0)=\frac12$ as $t_0\in U_m\setminus U_{m+1}$.
Let $h\in C_b(U)$. Let $V\subseteq U$ be an open neighbourhood of~$t_0$
such that $|h(t)-h(t_0)|<\eps$ for all $t\in V$.
Since $U_{m+1}$ is dense, $V\cap U_{m+1}\ne\emptyset$ and for each $t\in V\cap U_{m+1}$
we have $f_m(t)\geq\frac12+\frac14$ by construction; thus $|f_m(t)-f_m(t_0)|\geq\frac14$.
Suppose that $|f_m(s)-h(s)|\leq\eps$ for all $s\in U$. Then
\begin{equation*}
\textstyle{\frac14}\leq|f_m(t)-f_m(t_0)|\leq|f_m(t)-h(t)|+|h(t)-h(t_0)|+|h(t_0)-f_m(t_0)|<3\eps
\end{equation*}
wherefore $\eps>\frac{1}{12}$. As a result (\ref{equ:bigg_than_eps}) above holds for every
$\eps$ less than~$\frac{1}{12}$.

With the claim at hand we can now complete the proof as follows.
Let $H_n$, $n\in\NN$ be countably many mutually orthogonal infinite-dimensional
closed subspaces of~$H$ such that $H=\bigoplus_{n=1}^\infty H_n$.
Note that this induces a decomposition of $E$ as $E=\bigoplus_{n=1}^\infty E_n$
into an {\sl AW*}-direct sum, where $E_n$ is the {\sl AW*}-submodule of $E$
generated by the canonical isometric image of $H_n$ in~$E$.
Let $u_n\in B(H,H_n)$, $n\in\NN$ be unitaries and denote their extension
to isometric $B$-module isomorphisms from $E$ onto $E_n$ by the same symbol
so that $u_n(\rho(t))=u_n(\rho)(t)$ for all $\rho\in E$, $t\in X$. Put
$\eta_n=u_n(\zeta_n)$ for each~$n$. Then $\{\eta_n\mid n\in\NN\}$ is an orthonormal
family in the {\sl AW*}-module $E$ over~$B=C(X)$. By Lemma~\ref{lem:element-in-mke},
$a=\sum_{n=1}^\infty\eta_n\otimes\xi_n\in I(A)=\Mloc\mloc$, where we use the same
identification of $C(X,B(H)_\sigma)$ and $\mathcal L_B(E)$ as in
Proposition~\ref{prop:kb-in-mloc}. We shall show that $a\notin\mloc$.

Otherwise, by Corollary~\ref{cor:locmult-cxbh}, for given $\eps>0$
there are an open dense subset $U\subseteq X$ and $g\in C_b(U,B(H)_\beta)$
such that $\|a-g\|<\eps$. Denoting by $\tilde g$ the unique extension of $g$
to a function in $C(X,B(H)_\sigma)$ we find that
\[
\|a(t)-g(t)\|=\|a(t)-\tilde g(t)\|\leq\sup_{s\in X}\|a(s)-\tilde g(s)\|=\|a-g\|<\eps
\qquad(t\in U).
\]
Therefore,
\[
\|\eta_n(t)-g(t)(\xi_n)\|=\|a(t)(\xi_n)-g(t)(\xi_n)\|<\eps\qquad(t\in U).
\]
It follows that
\begin{equation*}
\begin{split}
\|\zeta_n(t)-u^*_np_ng(t)(\xi_n)\| &=\|u_n(\zeta_n(t))-p_ng(t)(\xi_n)\|\\
&\leq\|\eta_n(t)-g(t)(\xi_n)\|<\eps\qquad(t\in U),
\end{split}
\end{equation*}
where $p_n$ denotes the projection from $H$ onto~$H_n$.
Letting $g_n(t)=u_n^*p_ng(t)(\xi_n)$, $t\in U$, $n\in\NN$ we obtain a function
$g_n\in C_b(U,H)$ for each~$n$, since $g\in C_b(U,B(H)_\beta)$, such that
$\|\zeta_n-g_n\|<\eps$ for every $n\in\NN$. However, by the above claim, this is
impossible if $\eps<\frac{1}{12}$ and therefore $a\notin\mloc$.
\end{proof}

\begin{corollary}\label{cor:type-one-mloc}
There exists a unital type~{\rm I} \C*\/ $A$ such that $\mloc\ne\Mloc\mloc$.
\end{corollary}

\begin{proof}
Let $X$ be the spectrum of $\Mloc{C[0,1]}$; then $X$ has the properties required
in Theorem~\ref{thm:mloc-not-ia}, see the introductory remarks to this section.
Let $H=\ell^2$ and $C$ the unitisation of $K(H)$ inside $B(H)$. Then $A=C(X)\otimes C$
is a unital type~I \C*. Since $C(X)\otimes K(H)$ is an essential ideal in
both $A$ and $C(X)\otimes B(H)$, their local multiplier algebras agree with each other.
Thus the statement follows from Theorem~\ref{thm:mloc-not-ia}.
\end{proof}

\begin{remarks}
1.\ The proof of Theorem~\ref{thm:mloc-not-ia} uses some refinements of an argument
by Hamana to show statement~(2.5) on page~291 in~\cite{Hametc}.

\smallskip\noindent
2.\ One can add the condition of separability in Corollary~\ref{cor:type-one-mloc}
(and thus recover the example from~\cite{ArgFar2}, which was obtained
by those authors independently by different methods). Let $Y=[0,1]$. Using the
same reduction as in Corollary~\ref{cor:type-one-mloc} we can focus on the unital algebra
$A_0=C(Y)\otimes B(H)$, which is more in line with our previous discussion.

Letting $B=C(X)=\Mloc{C(Y)}$ as before we have
\begin{equation*}
\begin{split}
A=C(X)\otimes B(H)&=\dirlim\,C_b(V)\otimes B(H)\\
                  &=\dirlim\,\bigl(C(\beta V)\otimes B(H)\bigr)
                  \subseteq\dirlim\,C_b(V,B(H)_\beta)=\Mloc{A_0}
\end{split}
\end{equation*}
by Corollary~\ref{cor:locmult-cxbh} applied to~$A_0$, where the direct limit is taken over all
dense open subsets $V$ of~$Y$. Therefore,
$A_0\subseteq A\subseteq\Mloc{A_0}\subseteq I(A_0)$.
Applying Proposition~\ref{prop:invariant-injective} to $A_0$ we find that
$I(A_0)=I(A)$; hence the argument will be complete if we show that
$\Mloc{A_0}\subseteq\mloc$, since we already know that
$\mloc\ne I(A)=I(A_0)=\Mloc{\Mloc{A_0}}$, by \cite[Theorem~2.8]{Som2}.

This last assertion follows directly from the definition of~$\mloc$.
Recall that $X=\Invlim{\mathfrak D_Y}\beta V$, the projective limit along the filter
$\mathfrak D_Y$ of dense open subsets of~$Y$ and denote by
$\pi\colon X\to Y$ the canonical projection. For $V\in\mathfrak D_Y$, we have
$\pi^{-1}(V)\in\mathfrak D$ and since the embedding
$C_b(V,B(H)_\beta)\hookrightarrow C_b(\pi^{-1}(V),B(H)_\beta)\hookrightarrow\mloc$
is compatible with the connecting homomorphisms, we infer that
$\Mloc{A_0}\subseteq\mloc$.

The above observation answers Somerset's question on page~324 of~\cite{Som2}
in the negative.
\end{remarks}

We proceed to describe the maximal \C* of quotients of~$A$. Let $F=C(X,H)$, which is a Hilbert
$B$-submodule of $E=C(X,H_w)$.

\begin{proposition}\label{prop:closed-right-charac}
For a closed right ideal\/ $I$ of\/ $A$, let\/ $Z_I$ denote the closed
$B$-submodule of\/ $F$ generated by
\[
\{\rho\in F\mid\rho=a\xi\text{ for some }\,a\in I\text{ and some }\,\xi\in H\},
\]
where\/ $a\xi$ is the function\/ $(a\xi)(t)=a(t)\xi$ for every\/ $t\in X$.
For a closed $B$-submodule\/ $Z$ of\/ $F$, let\/ $I_Z=Z\otimes F$,
where, for\/ $\zeta\in Z$ and\/ $\rho\in F$,
$(\zeta\otimes\rho)(t)=\zeta(t)\otimes\rho(t)\in K(H)$ for all\/ $t\in X$. Then
\begin{enumerate}
\item[(i)] $Z_I$ is essential if\/ $I$ is essential;
\item[(ii)] $I_Z$ is a closed right ideal of\/ $A$ which is essential if\/ $Z$ is essential;
\item[(iii)] $I_{Z_I}\subseteq I$ for each\/ $I\in\Icr$.
\end{enumerate}
\end{proposition}

\begin{proof}
(i) Suppose that $I\in\Icer$ and let $0\ne\rho\in F$.
Take a unit vector $\eta\in H$. Since $I$ is essential,
one can find $\zeta\in F$ and $\xi\in H$, $\|\xi\|=1$ such that
$0\ne(\rho\otimes\eta)(\zeta\otimes\xi)\in I$. For each $t\in X$,
\[
(\rho(t)\otimes\eta)(\zeta(t)\otimes\xi)\xi=(\zeta(t)\mid\eta)\,\rho(t)
\]
and thus, with $f(t)=(\zeta(t)\mid\eta)$, $t\in X$ we obtain
$f\rho\in Z_I\setminus\{0\}$ for some $f\in B$ as desired.

(ii) For $\zeta\in Z$, $\rho\in F$ and $x\in A$ we have
\[
(\zeta\otimes\rho)\,x=\zeta\otimes x^*\rho\in Z\otimes F.
\]
Since $\{\zeta\otimes\rho\mid\zeta\in Z,\;\rho\in F\}$ is total in $Z\otimes F$,
this shows that $I_Z$ is a closed right ideal of~$A$.

Suppose that $Z$ is essential in~$F$. For a non-zero $x\in A$, there is $\xi\in H$
such that $x\xi\ne0$. Since $Z$ is essential, there is $h\in C(X)$ such that
$0\ne(x\xi)h\in Z$. Noting that
\[
x(t)(h(t)\xi\otimes\xi)= h(t)x(t)\xi\otimes\xi\qquad(t\in X)
\]
we conclude that $x(h\xi\otimes\xi)\in Z\otimes F\setminus\{0\}$,
which completes this step of the proof.

(iii) It suffices to show that $y\xi\otimes\rho\in I$ for every
$y\in I$, $\xi\in H$ and $\rho\in F$, which follows immediately from
\[
(y\xi\otimes\rho)(t)=y(t)(\xi\otimes\rho(t))=y(\xi\otimes\rho)(t)\qquad(t\in X).
\]
\vskip-10pt
\end{proof}

The proof of the following description of the essential
$C(X)$-submodules of~$F$ is straightforward and hence omitted.

\begin{proposition}\label{prop:essensubmodule}
A closed submodule of\/ $F$ is essential if and only if it contains a submodule
of the form\/ $\sum_{\rho\in F}C_0(U_\rho)\rho$, where, for each\/ $\rho\in F$,
the set\/ $U_\rho$ is open and dense in~$X$.
\end{proposition}

We are now ready to obtain our characterisation of~$\qmax$.

\begin{theorem}\label{thm:qmax-charac}
We have
\[
\qmaxsb=\bigl\{y\in C(X,B(H)_\sigma)\;\big|\;
               yF+y^*F\subseteq\Alglim{U\in\mathcal D}\,C_b(U,H)\bigr\}.
\]
\end{theorem}

\begin{proof}
Assume first that $y\in\qmaxsb$. Then there is a closed essential
right ideal $I$ of $A$ such that $yI+y^*I\subseteq A$.
By Proposition~\ref{prop:closed-right-charac},
there is an essential $B$-submodule $Z$ of $F$ such that $Z\otimes F\subseteq I$,
and, moreover, $Z\otimes F$ is an essential right ideal of~$A$. It is straightforward to
see that $yZ\subseteq F$. By Proposition~\ref{prop:essensubmodule},
for each $\rho\in F$, there is a dense open
subset $U_\rho$ of $X$ such that $C_0(U_\rho)\rho\subseteq Z$.
Given $\rho\in F$ we thus get that $C(X)\,C_0(U_\rho)\rho\subseteq C(X,H)$. We claim
that $y\rho$ is continuous on $U_\rho$. Indeed, for $t\in U_\rho$,
there is a function $g\in C_0(U_\rho)$ such that $g=1$ on a
neighbourhood of $t$ and, since $g(y\rho)$ is continuous by
hypothesis, we conclude that $y\rho$ is continuous at~$t$. This shows
that $y\rho\in\Alglim{U\in\mathcal D}\,C_b(U,H)$ and thus
$yF\subseteq\Alglim{U\in\mathcal D}\,C_b(U,H)$. An analogous argument holds for
$y^*F$.

Conversely, let $y\in C(X, B(H)_{\sigma})$ be such that
$yF+y^*F\subseteq\Alglim{U\in\mathcal D}\,C_b(U,H)$. For each $\rho\in F$,
there is $U_\rho\in \mathcal D$ such that the restriction of
$y\rho$ to $U_\rho$ is continuous. Let $Z$ be the closure of the
$B$-submodule $\sum_{\rho\in F}C_0(U_\rho)\rho$ of~$F$. By
Proposition~\ref{prop:essensubmodule}, $Z$ is a closed essential
submodule of $F$ and, by Proposition~\ref{prop:closed-right-charac},
we get that $I=Z\otimes F$ is an essential right ideal of~$A$. 
For $\rho\in F$ and $g\in C_0(U_\rho)$ we find that
\[
y(g\rho)=g(y\rho)\in C_0(U_\rho)\,C_b(U_\rho,H)\subseteq C_0(U_\rho,H)
\]
so that $yZ\subseteq F$, and consequently, $yI=yZ\otimes F\subseteq A$.
A similar argument shows that $y^*J\subseteq A$ for some $J\in\Icer$;
therefore $y$ belongs to~$\qmaxsb$.
\end{proof}

It has as yet still to be seen whether the result in Theorem~\ref{thm:qmax-charac}
suffices to determine whether $\qmax=I(A)$ if $A=C(X)\otimes B(H)$ in
the case where the Stonean space $X$ contains a dense subset of first category.

\section{Some Open Questions}\label{sect:problems}

In this section we list some open problems that arise from our discussion in the
previous sections.

\begin{question}\label{quest:two}
Following on from the last remark in the previous section, the following question
is close at hand.

\smallskip
\textit{With the notation and the assumptions of Theorem~\ref{thm:mloc-not-ia} is\/
$\qmax$ equal to or different from~$I(A)$?}
\end{question}

\begin{question}\label{quest:three}
Suppose that $\qmax\ne I(A)$ in the last problem; then the next question has a negative
answer, since, for $A=C(X)\otimes B(H)$, $\qmax$ contains all abelian projections
in $I(A)$ and thus would have to agree with $I(A)$ if it is an \AW*.

\smallskip
\textit{Is\/ $\qmax$ an \AW* for every \C*\/ $A$?}
\end{question}

\begin{question}\label{quest:four}
While we believe that the general answer to Question~\ref{quest:three} is negative,
the next problem should definitely have a positive answer. The discussion in
Section~\ref{sect:awstars} however shows that the present methods are not sufficient.

\smallskip
\textit{Does\/ $\qmax=A$ hold for every \AW*\/ $A$?}
\end{question}

\begin{question}\label{quest:five}
It was shown in~\cite{FrankPaulsen} that iterating the process of forming
the local multiplier algebra does not lead out of the injective envelope.
It is less clear what the relation to the maximal \C* of quotients may be.

\smallskip
\textit{Does\/ $M_{\text{\rm loc}}^{(k)}(A)\subseteq\qmax$ hold for every
\C*\/ $A$?}

\smallskip\noindent
Once again, we expect the answer to be `no'.
\end{question}

\begin{question}
It is known that $\Qmaxs\qmaxs=\qmaxs$ for every \C* $A$ from general algebraic results.
On the other hand, we know of examples in which $\Mloc\mloc\ne\mloc$. Hence we ask

\smallskip
\textit{Is it true that\/ $\Qmax\qmax=\qmax$ for every \C*\/ $A$?}
\end{question}

\end{document}